\documentclass[11pt, oneside]{amsart}
\usepackage[text={6.5in,9in},centering,letterpaper]{geometry}
\usepackage{graphicx}
\usepackage{amsfonts}
\usepackage[dvipsnames]{xcolor}
\usepackage{subcaption}
\usepackage{amssymb}
\usepackage{amsmath}
\usepackage{amscd}
\usepackage{pdfpages}
\usepackage{fancyhdr}
\usepackage{setspace}
\usepackage[all]{xy}
\usepackage{verbatim}
\usepackage{enumerate}
\usepackage[colorlinks=true, urlcolor=NavyBlue, linkcolor=NavyBlue, citecolor=NavyBlue,backref=page]{hyperref}
\usepackage{tikz}
\usetikzlibrary{  cd,  calc, positioning,  fit, arrows,  decorations.pathreplacing, decorations.markings,  shapes.geometric, backgrounds,bending}
\usepackage{tikzsymbols}
\usepackage{wrapfig}

\renewcommand*{\backref}[1]{}
\renewcommand*{\backrefalt}[4]
{%
    \ifcase #1 (Not cited.)%
        \or        (Cited on page~#2.)% cite of one page
        \else      (Cited on pages~#2.)% cite of several pages 
    \fi
}

\makeatletter

\renewcommand{\p@subfigure}{\thefigure.} % default: "\thefigure" (without the ".")
\makeatother

\newtheorem{theorem}{Theorem}[section]
\newtheorem{proposition}[theorem]{Proposition}
\newtheorem{lemma}[theorem]{Lemma}

\newtheorem{corollary}[theorem]{Corollary}

\makeatletter
\newtheorem*{rep@theorem}{\rep@title}
\newcommand{\newreptheorem}[2]{%
\newenvironment{rep#1}[1]{%
 \def\rep@title{#2 \ref{##1}}%
 \begin{rep@theorem}}%
 {\end{rep@theorem}}}
\makeatother

\newreptheorem{theorem}{Theorem}
\newreptheorem{lemma}{Lemma}
\newreptheorem{question}{Question}
\newreptheorem{corollary}{Corollary}
\newreptheorem{proposition}{Proposition}

\theoremstyle{definition}

\newtheorem{remark}[theorem]{Remark}

\newcommand{\Int}{\text{Int}}

\newcommand{\B}{\mathbb{B}}
\newcommand{\D}{\mathbb{D}}

\newcommand{\Z}{\mathbb{Z}}

\renewcommand{\S}{\mathbb{S}}

\newcommand{\Hh}{\mathcal H}

\newcommand{\Qq}{\mathcal Q}

\newcommand{\Ss}{\mathcal S}
\newcommand{\Tt}{\mathcal T}

\newcommand{\GG}{\mathfrak{G}}

\newcommand{\Sym}{\textsc{Sym}}

\newcommand{\zero}{\hbox{--\,}}

\makeatletter
\def\@seccntformat#1{%
  \protect\textup{\protect\@secnumfont
    \ifnum\pdfstrcmp{subsection}{#1}=0 \bfseries\fi% subsection # in \bfseries
    \csname the#1\endcsname
    \protect\@secnumpunct
  }%
}  
\makeatother

\begin{document}

\rhead{\thepage}
\lhead{\author}
\thispagestyle{empty}

%\tableofcontents
%\listoffigures

\raggedbottom
\pagenumbering{arabic}
\setcounter{section}{0}

%%%%%%%%%%%%%%%%%%%%%%%%%%%%%%%%%%%%%%%%%%%%%%%%%%%%%%%%
%%%%%%%%%%%%%%%%%%%%%%%%%%%%%%%%%%%%%%%%%%%%%%%%%%%%%%%%
%%%%%%%%%%%%%%%%%%%%%%%%%%%%%%%%%%%%%%%%%%%%%%%%%%%%%%%%

\title{A census of alternating tangles}
%\date{\today}

\author{Brendan D. McKay}
\address{School of Computing, Australian National University, Canberra, ACT, 2601, Australia}
\email{Brendan.McKay@anu.edu.au}
\urladdr{https://users.cecs.anu.edu.au/~bdm/}

\author{Jeffrey Meier}
\address{Department of Mathematics, Western Washington University, Bellingham, WA, 98229, USA}
\email{jeffrey.meier@wwu.edu}
\urladdr{http://jeffreymeier.org}

\begin{abstract}
    We give a census of the roughly 10 billion prime, connected, alternating tangles with at most 14 crossings.
    Along the way, we prove analogs of the Tait conjectures for alternating tangles and observe that aspects of the geometry of alternating tangles are reflected in their diagrams, as is the case with links.
\end{abstract}

\maketitle

%%%%%%%%%%%%%%%%%%%%%%%%%%%%%%%%%%%%%%%%%%%%%%%%%%%%%%%%
%%%%%%%%%%%%%%%%%%%%%%%%%%%%%%%%%%%%%%%%%%%%%%%%%%%%%%%%
\section{Introduction}
\label{sec:intro}
%%%%%%%%%%%%%%%%%%%%%%%%%%%%%%%%%%%%%%%%%%%%%%%%%%%%%%%%
%%%%%%%%%%%%%%%%%%%%%%%%%%%%%%%%%%%%%%%%%%%%%%%%%%%%%%%%

The purpose of this paper is to present a census of prime, connected, alternating tangles modulo isotopy rel. boundary and symmetry, in terms of their crossing number $c$.
Here, a tangle is a neatly-embedded 1--manifold $T$ in the $3$--ball $\B^3$; we allow any number $n\geq 2$ of arc-components, as well as any number $k\geq 0$ of circle-components.
The qualification that $T$ be prime precludes local knotting and split components, and the qualification that $T$ be connected means no diagram for $T$ can be separated by an arc in the disk; see Section~\ref{sec:anatomy} for precise definitions.

\begin{theorem}
\label{thm:counts}
    There are 10,070,064,432 prime, connected, alternating tangles with crossing number at most 14, up to isotopy rel. boundary and symmetry.
    Of these, 87,227,472 are diameter-prime, while 9,982,836,960 are not.
    The breakdown of these counts by crossing number is shown in Table~\ref{table:totals}.
\end{theorem}

\begin{table}[ht]
\setlength{\tabcolsep}{2.2pt}
\begin{tabular}{r|r|r|r|r|r|r|r|r|r|r|r|r|r|r}
$c$ & 1 & 2 & 3 & 4 & 5 & 6 & 7 & 8 & 9 & 10 & 11 & 12 & 13 & 14 \\
\hline
$N_1$ & 1 & 1 & 3 & 9 & 29 & 112 & 464 & 2200 & 11187 & 60819 & 344469 & 2013094 & 12009756 & 72785328 \\
$N_2$ & 0 & 1 & 3 & 16 & 88 & 588 & 4133 & 31025 & 239730 & 1901055 & 15350700 & 125903651 & 1046175570 & 8793230400 \\
\hline
Tot. & 1 & 2 & 6 & 25 & 117 & 700 & 4597 & 33225 & 250917 & 1961874 & 15695169 & 127916745 & 1058185326 & 8866015728
\end{tabular}
\caption{Counts of prime, connected, alternating tangles by crossing number~$c$. Indicated are the numbers of diameter-prime ($N_1$) and diameter-composite ($N_2$) tangles for each value of $c$.}
\label{table:totals}
\end{table}

Table~\ref{tab:counts} gives a breakdown of these counts (as well as counts for some other values of $c$ and $n$) in terms of the number of arc-components $n$ and circle-components $k$, and files containing the corresponding enumerations are publicly available in the Zenodo archive~\cite{zenodo}.
The algorithms used in generating these datasets are now implemented as a filter for the program \texttt{plantri}~\cite{plantri}; a complete description of this implementation is given in Section~\ref{sec:implementation}.

Our approach is to reduce the problem to the enumeration of quadrangulations of $\D^2$ modulo flype-equivalence and isomorphism.
Disk-quadrangulations are efficiently enumerated using the computer program \texttt{plantri}.
The main theoretical contribution is the following analog, for tangles and their diagrams, of the resolution of the Tait conjectures for knots and links; $c(T)$ denotes the crossing number of $T$ (see Section~\ref{sec:anatomy}).

\begin{theorem}
\label{thm:Tait_tangles}
	Let $D$ be a reduced, alternating diagram for a prime tangle $T$.
	Let $D'$ be any diagram for $T$.
	\begin{enumerate}
		\item If $D'$ is reduced and alternating, then $D$ and $D'$ are flype-equivalent modulo dihedral symmetry.
		\item $c(D')\geq c(D)$, with equality if and only if $D'$ is reduced and alternating.
	\end{enumerate}
\end{theorem}

We also observe that alternating tangles, by analogy with alternating links~\cite{menasco}, are prime, diameter-prime, and connected if and only if their alternating diagrams are; see Section~\ref{sec:geometry}.
Together, these results reduce the enumeration of prime (resp., diameter-prime), connected, alternating tangles to the enumeration of their alternating diagrams (Corollary~\ref{coro:tangle_to_diag}).

We consider tangles up to isotopy rel. boundary and symmetry; see Section~\ref{sec:anatomy} for details.
For prime, connected, alternating tangles, this reduces to the study of the underlying tangle projections modulo flype-equivalence and isomorphism (Lemma~\ref{lem:graph_to_diagram}).
Such tangle-graphs are in bijective correspondence with disk-quadrangulations, and this correspondence persists when working modulo flype-equivalence and isomorphism.
Thus, the enumeration of prime, connected, alternating tangles modulo isotopy rel. boundary and symmetry is accomplished via the enumeration of simple disk-quadrangulations modulo flype-equivalence and isomorphism, as codified in our final result.

\begin{repcorollary}{coro:tangle_to_quad}
    Every prime, connected, alternating tangle $T$ with $c(T)=c$ modulo isotopy rel. boundary and symmetry is represented uniquely by a simple quadrangulation $Q$ of $\D^2$ with $c$ faces modulo flype-equivalence and isomorphism.
    Moreover, $Q$ is diameter-prime if and only if $T$ is diameter-prime.
\end{repcorollary}

Our census divides prime, connected, alternating tangles based on whether they are diameter-prime or diameter-composite; see Section~\ref{sec:anatomy} for details.

%%%%%%%%%%%%%%%%%%%%%%%%%%%%%%%%%%%%%%%%%%%%%%%%%%%%%%%%
\subsection{Historical notes}
\label{subsec:history}
%%%%%%%%%%%%%%%%%%%%%%%%%%%%%%%%%%%%%%%%%%%%%%%%%%%%%%%%

The enumeration of prime, alternating knots and links via their diagrams has a long and rich history dating back to pioneering works of Tait, Kirkman, and Little from the end of the nineteenth century; we refer the reader to recent work by Burton~\cite{Bur_20_The-next-350-million-knots} and Thistlethwaite~\cite{Thi_25_The-enumeration-and-classification-of-prime-20-crossing} for a complete historical account.

A similar enumeration of alternating tangles was given by Bogdanov, Meshkov, Omelchenko, and Petrov~\cite{BogMesOme_12_Enumerating-the-k-tangle-projections}.
They enumerated prime, connected, alternating tangles through 12 crossings but omitted details regarding their method for accounting for flypes.
Moreover, they made implicit use of the ideas of Theorem~\ref{thm:Tait_tangles} without proof.
A completely different approach to the same problem was given by Jacobsen and Zinn-Justin~\cite{JacZin_02_The-combinatorics-of-alternating-tangles:-from}.
An analysis of the specific case of tangles with two arc-components was given by Sundberg and Thistlethwaite~\cite{SunThi_98_The-rate-of-growth-of-the-number-of-prime}, who proved a version of flype-equivalence in this setting.
These latter two works did not account for symmetry.

The enumeration of tangles up to isotopy was initiated by Kanenobu, Saito, and Satoh, who carried out their program up through seven crossings~\cite{KSS:seven}.

%%%%%%%%%%%%%%%%%%%%%%%%%%%%%%%%%%%%%%%%%%%%%%%%%%%%%%%%
\subsection*{Acknowledgements}
%%%%%%%%%%%%%%%%%%%%%%%%%%%%%%%%%%%%%%%%%%%%%%%%%%%%%%%%

The first author was supported by the Australian Research Council grant DP250101611.
The second author was supported by the U.S. National Science Foundation grant DMS-2405324.

%%%%%%%%%%%%%%%%%%%%%%%%%%%%%%%%%%%%%%%%%%%%%%%%%%%%%%%%%%%%
\subsection*{Statement of AI use}
%%%%%%%%%%%%%%%%%%%%%%%%%%%%%%%%%%%%%%%%%%%%%%%%%%%%%%%%%%%%

The authors used generative AI for help with proofreading, vector figure generation, coding, and background research.
The mathematical ideas presented in this article were conceived of and written by the authors without the aid of generative AI.
Generative AI was used to draw tangles, as in Figure~\ref{fig:tang2_6}, but was
not used in the generation of tangles.

%%%%%%%%%%%%%%%%%%%%%%%%%%%%%%%%%%%%%%%%%%%%%%%%%%%%%%%%
%%%%%%%%%%%%%%%%%%%%%%%%%%%%%%%%%%%%%%%%%%%%%%%%%%%%%%%%
\section{Tangle anatomy}
\label{sec:anatomy}
%%%%%%%%%%%%%%%%%%%%%%%%%%%%%%%%%%%%%%%%%%%%%%%%%%%%%%%%
%%%%%%%%%%%%%%%%%%%%%%%%%%%%%%%%%%%%%%%%%%%%%%%%%%%%%%%%

\begin{figure}[ht!]
    \centering
    \includegraphics[width=\textwidth]{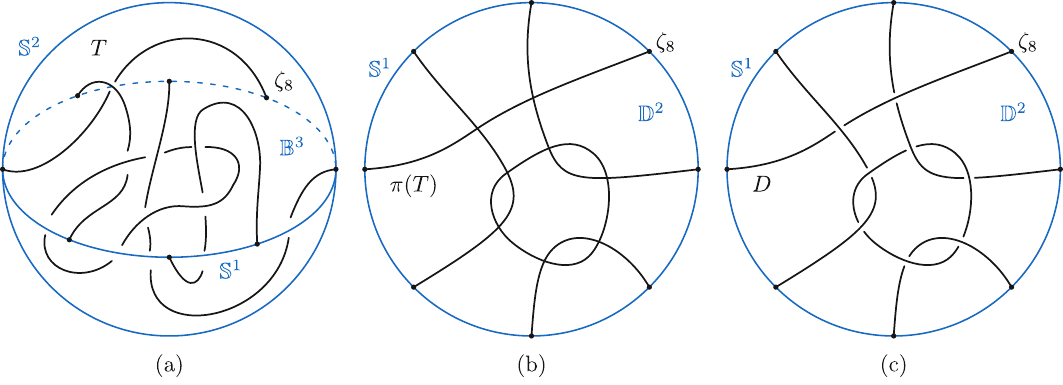}
    \caption{(a) A tangle $T$ in $\B^3$; (b) its projection $\Gamma = \pi(T)$ in $\D^2$; and (c) its diagram $D$.
    This tangle has $n=4$ arc-components and $k=1$ circle-component.
    The diagram $D$ is alternating and has $c(D)=8$ crossings.}
    \label{fig:tangle_ex}
\end{figure}

Let $(\B^3,T)$ be a tangle consisting of $n$ arc-components and $k$ circle-components for some $n\geq 2$ and $k\geq 0$.
We assume that the $2n$ points of $\partial T$ agree with the $2n^\text{th}$ roots of unity $\mathbf x_{2n} = \{\zeta_{2n}^i\}_{i=0}^{2n-1}$ on the boundary $\S^1$ of the equatorial unit disk $\D^2$.
After a small isotopy rel. boundary, we assume the projection $\pi\colon \B^3 \to \D^2$ is regular when restricted to $T$.
We call the graph $\Gamma = \pi(T)$ the \emph{projection} of $T$; this is a graph with $2n$ vertices of degree one (\emph{leaves}) at $\mathbf x_{2n}$ and a non-negative number of vertices of degree four (\emph{crossings}) in $\Int(\D^2)$.
We let $D$ denote the \emph{diagram} obtained from $\pi(T)$ in the usual way by recording at each crossing the relative height of the two strands of $T$ whose projections intersect at the crossing.
See Figure~\ref{fig:tangle_ex} for an example with $n=4$ and $k=1$.
A \emph{component} of $D$ is the portion of $D$ corresponding to a component of $T$.

Let $c(D)$ denote the number of crossings in a tangle diagram $D$.
Let $c(T)$ denote the \emph{crossing number} of a tangle $T$, which is defined as
$$c(T) = \min\{c(D)\,|\, \text{$D$ is a diagram for $T$}\}.$$

A tangle $T$ is \emph{prime} if any embedded 2--sphere in $\B^3$ intersecting $T$ in exactly two points bounds a 3--ball in $\B^3$ intersecting $T$ in an unknotted arc; otherwise, $T$ is \emph{composite}.
A tangle $T$ is \emph{split} if there exists an embedded 2--sphere $S\subset\B^3$ disjoint from $T$ that separates components of $T$; otherwise, $T$ is \emph{non-split}.

\begin{remark}
\label{rmk:n=1}
    If $T$ is a tangle with a single arc-component (so, $n=1$), then $\S^2$ intersects $T$ in two points.
    It follows that $T$ cannot be prime unless $T$ is an unknotted arc.
    However, we do not consider the unknotted arc to be a prime tangle, and we will assume throughout that $n\geq 2$ and $c(T)\geq 1$. 
    We avoid calling the unknotted arc the ``trivial tangle,'' because this term is typically applied to any tangle homeomorphic to a crossingless tangle.
\end{remark}

A neatly-embedded 2--disk $F\subset\B^3$ will be called \emph{vertical} if $\partial F\subset\S^2$ has a single maximum and a single minimum, and these critical points lie on opposite sides of the equator $\S^1$.
A tangle $T$ will be called \emph{diameter-split} if there is a vertical disk $F$ disjoint from $T$ that separates components of $T$; otherwise, $T$ will be called \emph{diameter-connected}.
A tangle $T$ will be called \emph{diameter-prime} if every vertical disk $F$ intersecting $T$ in exactly one point cuts off a 3--ball in $\B^3$ intersecting $T$ in an unknotted arc; otherwise, $T$ will be called \emph{diameter-composite}.

A tangle $T$ will be called \emph{connected} if it is non-split and diameter-connected.

\begin{remark}
\label{rmk:tangle_notions}
    For any tangle $T$ with $c(T)\geq 1$ or $n\geq 3$, we have the following implications:
$$
\begin{tikzcd}
\text{$T$ is split} \arrow[r, Rightarrow] \arrow[d, Rightarrow]
  & \text{$T$ is composite} \\
\text{$T$ is diameter-split} \arrow[r, Rightarrow]
  & \text{$T$ is diameter-composite}
\end{tikzcd}
$$
For the left implication, the splitting sphere can be tubed to a point on the equator $\S^1\subset\S^2$ along an arc disjoint from $T$ to give a splitting disk.
For the top implication, the splitting sphere can be moved by ambient isotopy (a finger-move) to intersect $T$ in two points connected by an unknotted arc; this is a decomposing sphere.
For the bottom implication, the splitting disk can be moved by an ambient isotopy (a finger-move along the boundary, through vertical disks) to intersect $T$ in a single point that is connected to $\S^2$ by an unknotted arc; this is a decomposing disk.
The bottom implication fails when $n=2$ and $c=0$.

In particular, $T$ is connected if and only if $T$ is diameter-connected.
\end{remark}

%%%%%%%%%%%%%%%%%%%%%%%%%%%%%%%%%%%%%%%%%%%%%%%%%%%%%%%%
\subsection{Tangle symmetry and equivalence}
\label{sec:symmetry}
%%%%%%%%%%%%%%%%%%%%%%%%%%%%%%%%%%%%%%%%%%%%%%%%%%%%%%%%

Let $\GG_{2n} = \Sym(\B^3,\mathbf x_{2n})$ denote the group of orthogonal symmetries of $\B^3$ that fix $\mathbf x_{2n}$ set-wise.
This group is isomorphic to $\Sym(\S^1,\mathbf x_{2n})\times\Z_2$, where $\Sym(\S^1,\mathbf x_{2n})$ is the dihedral group of symmetries of the regular $2n$--gon, which we regard as having vertices $\mathbf x_{2n}$.
The group $\GG_{2n}$ has $8n$ elements and is generated by three transformations:
\begin{enumerate}
    \item the rotation $\rho$ through $\pi/n$ radians that is right-handed about the $z$--axis in $\B^3$;
    \item the flip $\varphi$ that rotates $\B^3$ through $\pi$ radians about a diameter contained in $\D^2$ and intersecting $\zeta_{2n}^0\in\mathbf x_{2n}$; and
    \item the mirroring $\mu$, which is reflection of $\B^3$ across $\D^2$.
\end{enumerate}
At the level of a tangle diagram $D$, the transformations $\rho$, $\varphi$, and $\mu$ have the effect (respectively) of rotating $D$ counterclockwise through $\pi/n$ radians, reflecting $D$ across the horizontal diameter of $\D^2$ before reversing each crossing, and reversing each crossing of $D$.
For any diagram $D$, we denote the \emph{mirror} diagram $\mu(D)$ by $\overline D$.

Two tangles are \emph{equivalent} if they are in the same orbit of the action of $\GG_{2n}$, up to isotopy rel. boundary.
If tangles are in the same orbit of the subgroup of $\GG_{2n}$ generated by $\rho$ and $\varphi$, up to isotopy rel. boundary, we say they are equivalent modulo \emph{dihedral symmetry}.
Equivalent tangles have diagrams that are related by Reidemeister moves and symmetry.

%%%%%%%%%%%%%%%%%%%%%%%%%%%%%%%%%%%%%%%%%%%%%%%%%%%%%%%%
\subsection{Connectedness in tangle diagrams}
\label{sec:connected_diagrams}
%%%%%%%%%%%%%%%%%%%%%%%%%%%%%%%%%%%%%%%%%%%%%%%%%%%%%%%%

A tangle diagram $D$ is \emph{split} if there is a simple closed curve $\gamma\subset\D^2$ disjoint from $\Gamma$ that separates components of $\Gamma$; otherwise, $D$ is \emph{non-split}.
A tangle diagram $D$ is \emph{prime} if any simple closed curve $\gamma\subset\D^2$ intersecting $\Gamma$ in exactly two (non-crossing) points bounds a disk in $\D^2$ intersecting $D$ in a diagram for an unknotted arc; otherwise, $D$ is \emph{composite}.

A tangle diagram $D$ is \emph{diameter-split} if there is a neatly-embedded arc $\gamma\subset\D^2$ disjoint from $\Gamma$ that separates components of $\Gamma$; otherwise, $D$ is \emph{diameter-connected}.
A tangle diagram $D$ is \emph{diameter-prime} if any neatly-embedded arc $\gamma\subset\D^2$ intersecting $\Gamma$ in exactly one (non-crossing) point bounds a disk to one side intersecting $D$ in a diagram for an unknotted arc; otherwise, $D$ is \emph{diameter-composite}.

A tangle diagram $D$ is \emph{connected} if its underlying projection $\Gamma$ is connected (as a graph); otherwise, $D$ is \emph{disconnected}.
Note that a tangle diagram $D$ is connected if and only if it is diameter-connected and non-split.

\begin{remark}
\label{rmk:diagram_notions}
    The diagrammatic notions defined immediately above for tangle diagrams correlate with the corresponding geometric notions for tangles defined earlier in this section to the following extent.
    If $D$ is a diagram for $T$ and $D$ is split, composite, diameter-split, diameter-composite, or disconnected, then $T$ is split, composite, diameter-split, diameter-composite, or disconnected, respectively.
    To see this, simply note that a splitting/decomposing circle in $D$ can be capped off with disks above and below $\D^2$ to give a splitting/decomposing sphere for $T$, and a splitting/decomposing arc for $D$ can be vertically extruded to give a splitting/decomposing disk for $T$.
    In Section~\ref{sec:geometry}, we will show that these geometric properties of $T$ are apparent in $D$ whenever $D$ is alternating.

    Note that $D$ being split does not necessarily imply $D$ is diameter-split, in contrast to the left implication in Remark~\ref{rmk:tangle_notions}.
    However, $D$ is diameter-composite whenever $D$ is diameter-split and either $n\geq 3$ or $c(T)\geq 1$, and $D$ is composite whenever $D$ is split.
\end{remark}

A crossing $x$ of a tangle diagram $D$ is \emph{nugatory} if there is a simple closed curve in $\D^2$ intersecting $D$ in exactly one point: the crossing $x$.
We say that $D$ is \emph{reduced} if it has no nugatory crossings.

A tangle diagram $D$ is \emph{strongly prime} if any simple closed curve $\gamma\subset\D^2$ intersecting $D$ in exactly two (non-crossing) points bounds a disk intersecting $D$ in an embedded arc.

%%%%%%%%%%%%%%%%%%%%%%%%%%%%%%%%%%%%%%%%%%%%%%%%%%%%%%%%
\subsection{Alternating tangle diagrams and flyping}
\label{sec:flyping}
%%%%%%%%%%%%%%%%%%%%%%%%%%%%%%%%%%%%%%%%%%%%%%%%%%%%%%%%

A modicum of care must be taken when defining what it means for a tangle to be alternating if disconnected tangles are to be considered.
For our purposes, we will require that all arc-components of the diagram are alternating in a consistent way.
We will call a tangle diagram \emph{alternating} if the following two conditions hold:
\begin{enumerate}
    \item As one travels along any component of the diagram, one encounters crossings alternately as the over-strand or under-strand, or one encounters no crossings at all; and
    \item If the first crossing encountered when traveling along an arc-component starting at the endpoint $\zeta_{2n}^i$ for some $0\leq i<2n$ is encountered as the over-strand, then the same is true whenever a first crossing is encountered when starting such a journey from the endpoint $\zeta_{2n}^j$ for any $j\equiv i\pmod 2$.
\end{enumerate}
By a checkerboard argument, condition (2) implies that any first crossing encountered when traveling along an arc-component starting at the endpoint $\zeta_{2n}^j$ for any $j\not\equiv i\pmod 2$ is encountered as the under-strand.
The diagram shown in Figure~\ref{fig:tangle_ex}(c) is alternating.
A tangle is \emph{alternating} if it has an alternating diagram.

Let $\Ss\subset\S^3$ denote the one-point compactification of the $xy$--plane, and write $\Ss = \D^2\cup_{\S^1}\D^2_\infty$, where $\D^2_\infty$ is the complementary disk containing $\infty$.
Note that with the present definition, a tangle diagram $D\subset\D^2$ is alternating if and only if every link diagram $\widehat D$ obtained by connecting the endpoints of $D$ with $n$ arcs in $\D^2_\infty$ that do not cross is an alternating link diagram.
Under this definition, a non-split (but possibly diameter-split) projection $\Gamma$ can be made alternating in exactly two ways, rather than $2^m$ ways, where $m$ is the number of diameter-split components.
(Recall that we assume throughout that $\Gamma$ has at least one crossing.)

\begin{figure}[ht!]
    \centering
    \includegraphics[width=.8\textwidth]{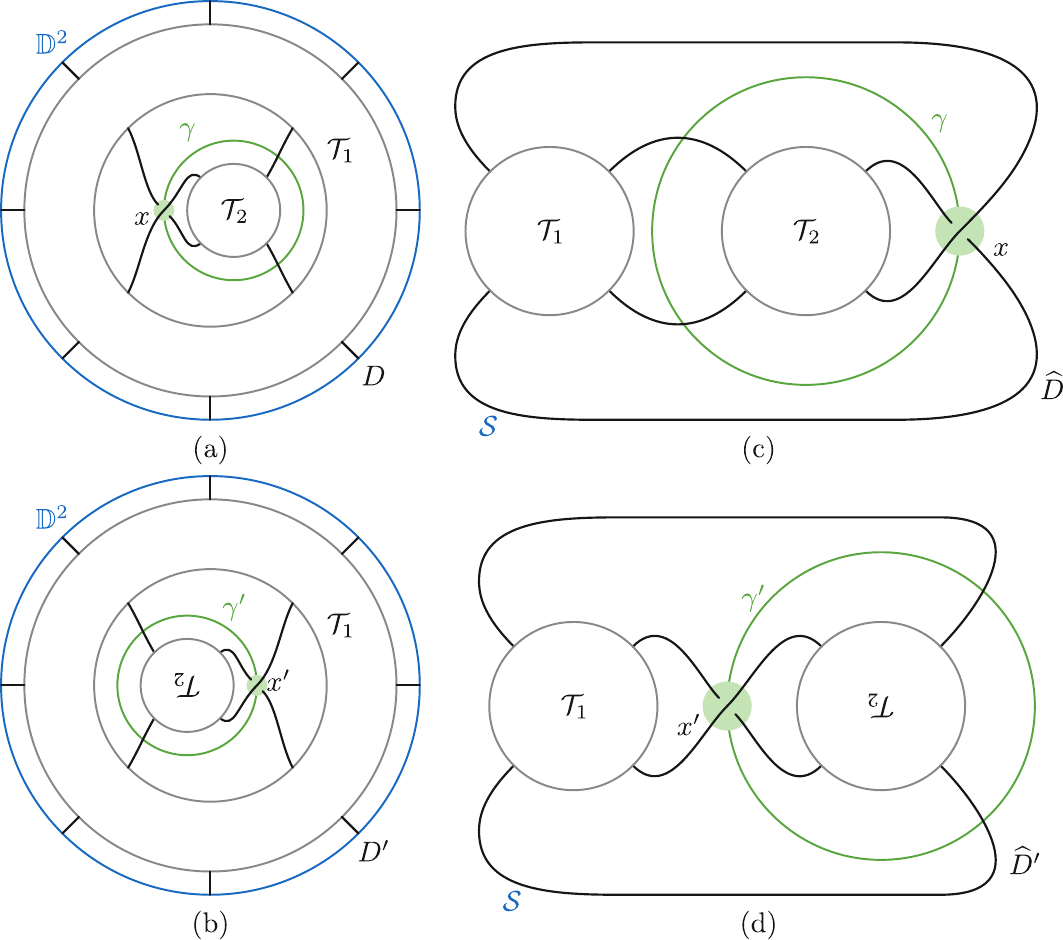}
    \caption{(a)--(b) A flype of a tangle diagram $D$ (along $\gamma$) commuting the crossing $x$ in $D$ to the crossing $x'$ in the resulting diagram $D'$; and (c)--(d) a flype of a link diagram $\widehat D$ (along $\gamma$) commuting the crossing $x$ in $\widehat D$ to the crossing $x'$ in the resulting diagram $\widehat D'$.}
    \label{fig:flype_def}
\end{figure}

Given a tangle diagram $D$, a simple closed curve $\gamma\subset\D^2$ is a \emph{flyping circle} for $D$ if it intersects $D$ in three points, exactly one of which is a crossing $x$, which it intersects transversely, dividing the four half-edges into two adjacent pairs; see Figure~\ref{fig:flype_def}(a).
A \emph{flype} (along $\gamma$) is the geometric operation whereby the portion of $D$ inside $\gamma$ is rotated through $\pi$ radians in such a way that the crossing $x$ is commuted across the rest of the region bounded by $\gamma$, as indicated in Figures~\ref{fig:flype_def}(a) and~\ref{fig:flype_def}(b), where the portion of $D$ labeled $\Tt_2$ has been flipped upside down.

\begin{remark}
\label{rmk:flip_flype}
    It is also possible to flip over the outer, annular region $\Tt_1$.
    In this case, we assume the flip preserves $\mathbf x_{2n}$ set-wise, so it fixes antipodal points $y$ and $y'$, which are either both in $\mathbf x_{2n}$ or both midpoints of points in $\mathbf x_{2n}$. 
    After a planar isotopy rel. boundary of $D$, we can assume the fixed-set of this flip lies on a diameter of $\D^2$ connecting $y$ and $y'$.
    The result of following the flip of $\Tt_1$ with the global flip about this diameter is precisely the flype that flips over the inner region $\Tt_2$.
    Thus, the two different flypes along $\gamma$ are related by a dihedral symmetry.
    Since we are only concerned with studying tangle diagrams modulo symmetry, the choice of region does not affect the end result.
\end{remark}

Two diagrams $D$ and $D'$ are \emph{flype-equivalent} if they are related by a sequence of flypes.

%%%%%%%%%%%%%%%%%%%%%%%%%%%%%%%%%%%%%%%%%%%%%%%%%%%%%%%%
\subsection{Corona-closures}
\label{sec:corona-closure}
%%%%%%%%%%%%%%%%%%%%%%%%%%%%%%%%%%%%%%%%%%%%%%%%%%%%%%%%

\begin{figure}[ht!]
	\centering
	\includegraphics[width = \textwidth]{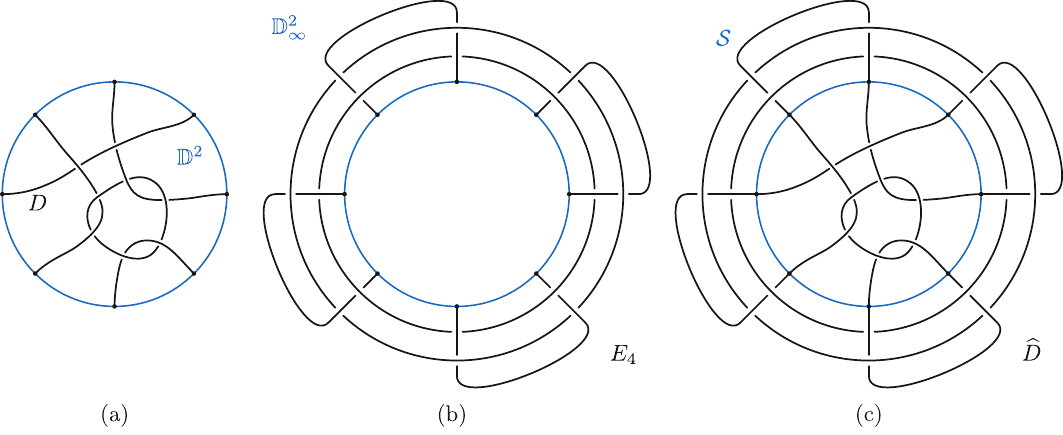}
	\caption{(a) A tangle diagram $D$ in $\D^2$; (b) a corona diagram $E_4$ in $\D^2_\infty$; and (c) the diagram $\widehat D$ in $\Ss$ for the corona-closure of $D$.}
	\label{fig:corona}
\end{figure}

We now introduce a bespoke tangle that will play a central role in many of the proofs that follow.
Let $D_n\subset\D^2$ denote the alternating tangle diagram with $n+2$ components built as follows:
\begin{enumerate}
    \item Take $n$ pairwise disjoint arc-components $a_1,\ldots, a_n$ such that $a_i$ goes from $\zeta_{2n}^{2i-2}$ to $\zeta_{2n}^{2i-1}$.
    \item Add two circle-components $c_1$ and $c_2$ that are nested, concentric copies of $\S^1$, each of which crosses each $a_i$ exactly twice.
    \item Ensure that the diagram is alternating such that $a_1$ is the over-strand of the crossing nearest to $\zeta_{2n}^0$.
\end{enumerate}
Finally, let $E_n\subset\D^2_\infty$ denote the result of reflecting $D_n$ across $\S^1$.
We call the diagram $E_n$ a \emph{corona diagram}; see Figure~\ref{fig:corona}(b) for the case of $n=4$.

Given a tangle diagram $D$ with $n$ arc-components, consider the link diagram in $\Ss$ given by $\widehat D = D\cup E_n$.
We call $\widehat D$ the \emph{corona-closure} of $D$; see Figure~\ref{fig:corona} for an example of the process of converting a tangle diagram $D$ to its corona-closure.
Note that $E_n$ is visibly prime: No pair of polygonal regions of $\D^2_\infty$ separated by the diagram share two distinct boundary edges.

The following key lemma relates properties of a tangle diagram $D$ to corresponding properties of its corona-closure $\widehat D$, which we now define.
A link diagram $\widehat D$ is \emph{prime} if any simple closed curve in $\Ss$ that intersects $\widehat D$ in two points (that are not crossings) bounds a disk intersecting $\widehat D$ in a diagram for an unknotted arc; otherwise, $\widehat D$ is \emph{composite}.
A link diagram $\widehat D$ is \emph{split} if there is a simple closed curve in $\Ss$ disjoint from $\widehat D$ and separating components of $\widehat D$.
A link diagram is \emph{connected} if it is non-split.

A crossing $x$ of a link diagram $\widehat D$ is \emph{nugatory} if there is a simple closed curve in $\Ss$ intersecting $\widehat D$ in exactly one point: the crossing $x$.
We say that $\widehat D$ is \emph{reduced} if it has no nugatory crossings.
A simple closed curve $\gamma\subset\Ss$ is a \emph{flyping circle} for a link diagram $\widehat D$ if it intersects $\widehat D$ in three points, exactly one of which is a crossing $x$, which it intersects transversely, dividing the four half-edges into two adjacent pairs; see Figure~\ref{fig:flype_def}(c).
A \emph{flype} (along $\gamma$) is the geometric operation whereby the portion of $\widehat D$ inside $\gamma$ is rotated through $\pi$ radians in such a way that the crossing $x$ is commuted across the rest of the region bounded by $\gamma$, as indicated in Figures~\ref{fig:flype_def}(c) and~\ref{fig:flype_def}(d), where the portion of $\widehat D$ labeled $\Tt_2$ has been flipped upside down.
It is also possible to flip over the outer region $\Tt_1$.

\begin{lemma}
\label{lem:corona}
    Let $D$ be a tangle diagram, and let $\widehat D$ be its corona-closure.
    Then,
    \begin{enumerate}
        \item $D$ is split, prime, or reduced if and only if $\widehat D$ is split, prime, or reduced, respectively.
        \item If $D$ is strongly prime, then any flyping circle for $\widehat D$ either intersects $\D^2_\infty$ as shown in Figure~\ref{fig:local}(f) or else is ambiently isotopic (through flyping circles) to a flyping circle contained in either $\D^2$ or $\D^2_\infty$.
        \item If $D$ is alternating, then either $\widehat D$ is alternating, or the corona-closure of the mirror diagram $\overline D$ is alternating.
    \end{enumerate}
\end{lemma}

\begin{proof}
    In part (1), the fact that $\widehat D$ is split, is composite, or has a nugatory crossing whenever $D$ is split, is composite, or has  a nugatory crossing, respectively, is immediate.
    
    For the first converse, assume $\gamma$ is a splitting circle for the split diagram $\widehat D$.
    Since $E$ is non-split, any component of $\D^2_\infty\cap\gamma$ is a boundary-parallel arc of the type shown in Figure~\ref{fig:local}(a).
    Any such arc can be moved by an ambient isotopy through splitting circles to lie in $\D^2$.
    After such modification, $\gamma$ is a splitting circle for $D$.

\begin{figure}[ht!]
    \centering
    \includegraphics[width = \textwidth]{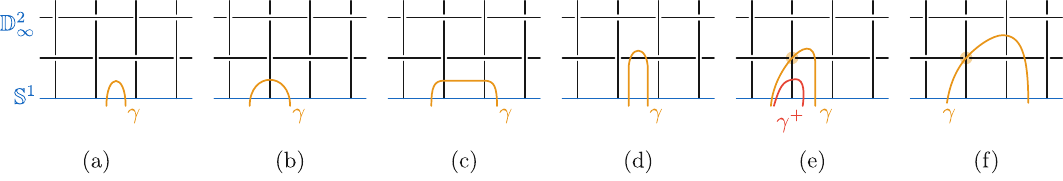}
    \caption{Local pictures for the intersection of a splitting, decomposing, reducing, or flyping circle $\gamma$ with $\D^2_\infty$ in the diagram for a corona-closure. The corona-closure caps are suppressed for generality.}
    \label{fig:local}
\end{figure}

    For the second converse, assume $\gamma$ is a decomposing circle for the composite diagram $\widehat D$.
    Since $E$ is prime, every component of $\D^2_\infty\cap\gamma$ has one of the forms shown in Figures~\ref{fig:local}(a) through~\ref{fig:local}(d).
    If the intersection is of the first three forms shown, then $\gamma$ can be moved by an ambient isotopy through decomposing circles to lie in $\D^2$ as a decomposing circle for $D$.
    If the intersection is of the fourth form, then $\gamma$ is disjoint from $D$, so $D$ must be split.
    But this means $D$ is composite, since the splitting circle is ambient isotopic to a decomposing circle.

    For the third converse, assume $x$ is a nugatory crossing for the diagram $\widehat D$, and let $\gamma$ be the reducing circle in $\Ss$ intersecting $\widehat D$ in exactly one point: the crossing $x$.
    Note that $x$ cannot lie in $\D^2_\infty$ since all such crossings involve distinct components.
    So, $x$ must be a crossing of $D$ and any portion of  $\D^2_\infty\cap\gamma$ must be of the form shown in Figure~\ref{fig:local}(a).
    It follows that $\gamma$ can be moved by ambient isotopy to lie in $\D^2$, so $x$ is a nugatory crossing for $D$.

    For part (2), let $\gamma$ be a flyping circle for $\widehat D$.
    If $\gamma\subset\D^2_\infty$ or $\gamma\subset\D^2$, then we are done.
    Otherwise, every component of $\D^2_\infty\cap\gamma$ has one of the six forms shown in Figure~\ref{fig:local}.
    Any occurrences of the first three types can be moved by an ambient isotopy through flyping circles to lie in $\D^2$.
    Any occurrence of the fourth type is impossible:
    By a (mod 2)--intersection argument, the crossing $x$ in $\D^2$ must be a nugatory crossing, but this contradicts the assumption that $D$ is strongly prime.
    So we can assume that $\D^2_\infty\cap\gamma$ is a single instance of the fifth or sixth form.
    In the latter case, we are done.
    In the former case, let $\gamma^+$ be a slight push-off of $\gamma$, as shown in Figure~\ref{fig:local}(e).
    Then, $\gamma^+$ is ambiently isotopic into $\D^2$ through circles intersecting $\widehat D$ in two points.
    Since $D$ is strongly prime, $\gamma$ is ambiently isotopic through flyping circles to a flyping circle contained in $\D^2_\infty$.

    Part (3) follows from condition (2) in the definition of alternating tangle diagrams and the comment immediately following it regarding a checkerboard argument.
\end{proof}

%%%%%%%%%%%%%%%%%%%%%%%%%%%%%%%%%%%%%%%%%%%%%%%%%%%%%%%%
%%%%%%%%%%%%%%%%%%%%%%%%%%%%%%%%%%%%%%%%%%%%%%%%%%%%%%%%
\section{Alternating tangle diagrams}
\label{sec:alternating_tangles}
%%%%%%%%%%%%%%%%%%%%%%%%%%%%%%%%%%%%%%%%%%%%%%%%%%%%%%%%
%%%%%%%%%%%%%%%%%%%%%%%%%%%%%%%%%%%%%%%%%%%%%%%%%%%%%%%%

In this section, we first recall three theorems regarding alternating links that were conjectured by P. G. Tait circa 1876 and were proved by various authors circa 1990.
We then prove that analogous theorems hold for alternating tangles.
The terms ``prime'' and ``strongly prime'' are not consistently applied in the literature, and some theorems were originally proved under the stronger hypothesis of strong primeness.
For a treatment of the classical theorems that is consistent with our approach and definitions, see~\cite[Chapters~4--5]{Lic_97_An-introduction-to-knot-theory}, especially Theorem~4.3 and Corollary~5.10.

The first theorem concerns crossing-minimal diagrams of alternating links.

\begin{theorem}[\cite{Kau_87_State-models-and-the-Jones-polynomial, Mur_87_Jones-polynomials-and-classical-conjectures, Thi_87_A-spanning-tree-expansion-of-the-Jones, Tur_87_A-simple-proof-of-the-Murasugi-and-Kauffman}]
\label{thm:tait_1}
	If $\widehat D$ is a reduced, alternating diagram for a link $L$, then no diagram for $L$ has fewer crossings than $\widehat D$.
	Moreover, any non-alternating, prime diagram for $L$ has more crossings than $\widehat D$.
\end{theorem}

The second theorem concerns flype-equivalence of alternating diagrams.
Recall that a link $L\subset\S^3$ is \emph{split} if there is an embedded 2--sphere that is disjoint from $L$ and separates components of $L$.
A link $L\subset\S^3$ is \emph{prime} if any embedded 2--sphere intersecting $L$ in exactly two points bounds a ball intersecting $L$ in an unknotted arc.
Note that the only prime, split link is the unlink of two components.

\begin{theorem}[\cite{MenThi_91_The-Tait-flyping-conjecture, MenThi_93_The-classification-of-alternating-links}]
\label{thm:tait_2}
	Any two reduced, alternating diagrams for a given prime link are related by a sequence of flypes.
\end{theorem}

If $L$ (and $\widehat D$) are oriented, then the \emph{writhe} of $\widehat D$ is the signed count of the crossings of $\widehat D$ based on the local orientations; see~\cite[Figure~1.11]{Lic_97_An-introduction-to-knot-theory}.
The count is not invariant across diagrams for $L$ in general.
This third theorem follows from Theorem~\ref{thm:tait_2}, since flyping preserves the writhe of a diagram, but was originally proved using the same techniques as Theorem~\ref{thm:tait_1}.

\begin{theorem}[\cite{Mur_87ii_Jones-polynomials-and-classical-conjectures, Thi_88_Kauffmans-polynomial-and-alternating-links}]
\label{thm:tait_3}
	Any two reduced, alternating diagrams for a given oriented, prime link have the same writhe.
\end{theorem}

%%%%%%%%%%%%%%%%%%%%%%%%%%%%%%%%%%%%%%%%%%%%%%%%%%%%%%%%
\subsection{Alternating tangle diagrams}
\label{subsec:tangle_diagrams}
%%%%%%%%%%%%%%%%%%%%%%%%%%%%%%%%%%%%%%%%%%%%%%%%%%%%%%%%

The theorems about alternating tangle diagrams proved in this section are analogs of, and follow from, the corresponding theorems about alternating link diagrams given above.
The general approach is to turn an alternating tangle diagram into an alternating link diagram by considering a corona-closure.

Our first theorem is the analog of Theorem~\ref{thm:tait_1} and concerns crossing-minimal diagrams for alternating tangles.

\begin{theorem}
\label{thm:alt_tang_1}
	If $D$ is a reduced, alternating diagram for a tangle $T$, then no diagram for $T$ has fewer crossings than $D$.
	Moreover, any non-alternating, prime diagram for $T$ has more crossings than $D$.
\end{theorem}

\begin{proof}
	Let $D$ be a reduced, alternating diagram for a given tangle $T$.
	Let $E$ be the corona diagram $E_n$, where $n$ is the number of arc-components of $D$.

	First, note that at least one of $D\cup E$ and $\overline D\cup E$ will be alternating, by Lemma~\ref{lem:corona}(3).
	Without loss of generality, assume $\widehat D = D\cup E$ is alternating.
	
	Suppose $\widehat D$ is a diagram for a link $L$.
	By Lemma~\ref{lem:corona}(1), since $D$ is reduced, so is $\widehat D$.
	By Theorem~\ref{thm:tait_1}, it follows that $\widehat D$ minimizes crossing number for $L$.
	
	If $D'$ were a diagram for $T$ with fewer crossings than $D$, then $\widehat D' = D'\cup E$ would be a diagram for $L$ with fewer crossings than $\widehat D$, a contradiction.
	Therefore, no diagram for $T$ has fewer crossings than~$D$.
	
	Now, suppose $D'$ is a prime, non-alternating diagram for $T$.
	Then $\widehat D'$ is a non-alternating diagram for $L$.
    By Lemma~\ref{lem:corona}(1), the fact that $D'$ is prime implies that $\widehat D'$ is prime.
	By Theorem~\ref{thm:tait_1}, $\widehat D'$ must have more crossings than $\widehat D$, so $D'$ must have more crossings than $D$.
\end{proof}

Our second theorem is the analog of Theorem~\ref{thm:tait_2} and concerns flype-equivalence of alternating tangle diagrams.

\begin{theorem}
\label{thm:alt_tang_2}
	Any two reduced, alternating diagrams for a given prime tangle are related by a sequence of flypes and dihedral symmetry.
\end{theorem}

\begin{proof}
	Let $T$ be a prime tangle, and let $D$ and $D'$ be reduced, alternating diagrams for $T$.
	
	As in the previous proof, by Lemma~\ref{lem:corona}(3), assume without loss of generality that $\widehat D = D\cup E$ is alternating.
	We claim that $\widehat D' = D'\cup E$ is also alternating.
    Suppose for a contradiction that it is not.
    
	By Theorem~\ref{thm:alt_tang_1}, $D$ and $D'$ have the same number of crossings, so the link diagrams $\widehat D'$ and $\widehat D$ have the same number of crossings.
    Since $T$ is prime, so are $D$ and $D'$.
    By Lemma~\ref{lem:corona}(1), this implies $\widehat D$ and $\widehat D'$ are prime.
	However, by assumption, $\widehat D'$ is non-alternating, so by Theorem~\ref{thm:tait_1}, it has more crossings than $\widehat D$, a contradiction.
	Therefore, we conclude that if $\widehat D$ is alternating, so is $\widehat D'$, as claimed.
	
	So, $\widehat D$ and $\widehat D'$ are alternating diagrams for the same link $L$.
	By Lemma~\ref{lem:corona}(1), since $D$ and $D'$ are reduced, the same is true for $\widehat D$ and $\widehat D'$.
	We claim that since $T$ is prime, so is $L$.
    To see this, suppose $L$ is composite.
	By~\cite[Theorem~4.4]{Lic_97_An-introduction-to-knot-theory} (see Theorem~\ref{thm:links_split}, below), this implies $\widehat D$ (and $\widehat D'$) are composite.
    By Lemma~\ref{lem:corona}(1), $D$ (and $D'$) are composite, and this implies $T$ is composite, by Remark~\ref{rmk:diagram_notions}.

    Because $\widehat D$ and $\widehat D'$ are reduced, alternating diagrams for the prime link $L$, by Theorem~\ref{thm:tait_2}, $\widehat D'$ is obtained from $\widehat D$ by a sequence of flypes.
    In fact, Menasco and Thistlethwaite prove a stronger statement than Theorem~\ref{thm:tait_2}:
    They prove that any orientation-preserving homeomorphism of pairs $f\colon(\S^3,L)\to(\S^3,L')$, where $L$ and $L'$ are oriented links having diagrams $\widehat D$ and $\widehat D'$, respectively, is pairwise-isotopic to a composition of flypes and flat homeomorphisms~\cite[Theorem~1]{MenThi_91_The-Tait-flyping-conjecture}.
    For our purposes, since $D$ and $D'$ are both diagrams for $T$, we can assume that $f$ restricts to the identity outside of $\B^3$.
    In particular, $f$ preserves (with labels) the two circle-components of $E$.
    It follows that we can assume that the sequence of flypes converting $\widehat D$ into $\widehat D'$ preserves (with labels) the two circle-components of $E$, as well.
    
    We will use this extra structure on $f$ at the end of the proof.
    For now, we prove the following claim, which will be applied to each flype in the sequence relating $\widehat D$ and $\widehat D'$.
    Note that $D$ is strongly prime since it is prime and crossing-minimal.

    \textbf{Claim:} If $\gamma$ is a flyping circle for $\widehat D = D\cup E$, with $D$ strongly prime, then the result of a flype along $\gamma$ is $\widehat D'' = D''\cup E$ where $D''$ is strongly prime and obtained from $D$ by (a) a flype, (b) a dihedral symmetry (including the identity), or (c) a flype and a dihedral symmetry.

    \begin{proof}[Proof of Claim]
    To prove this claim, first apply Lemma~\ref{lem:corona}(2) to ensure that either $\D^2_\infty\cap\gamma$ has the form shown in Figure~\ref{fig:local}(f) or else $\gamma$ lies entirely in either $\D^2$ or $\D^2_\infty$.
    If $\gamma\subset\D^2$, then the flype of $\widehat D$ along $\gamma$ restricts to a flype of $D$ along $\gamma$.
    If the flype flips over the inside region, take $D''$ to be the result of applying the flype to $D$.
    If the outside region is the one that is flipped, the diagram $E$ is flipped about a diameter of $\D^2_\infty$.
    The result is either $E$ or $\rho(E)$.
    In the former case, take $D''$ to be the result of the flype on $D$ that flips the outside, annular region.
    In the latter case, take $D''$ to be the result of this same flype, followed by an application of $\rho^{-1}$.
    
    If $\gamma\subset\D^2_\infty$, then $\gamma$ bounds a disk in $\D^2_\infty$ containing one or fewer crossings of $\widehat D$ (not counting $x$).
    In any case, the flype that flips over the inside region does not change the diagram $\widehat D$, so we can take $D''=D$, and the flype that flips over the outside region flips $\D^2$ along a diameter.
    After a potential rotation, the corona is restored, so we can take $D''$ to be the result of applying this total dihedral symmetry to $D$.

    Finally, we consider the case where $\D^2_\infty\cap\gamma$ has the form shown in Figure~\ref{fig:local}(f).
    In this case, $\D^2\cap\gamma$ is a splitting arc for $D$, so $\gamma$ bounds a disk intersecting $\widehat D$ in a single crossing (not counting $x$) and intersecting $D$ in a single arc.
    This arc is crossingless, since $D$ is strongly prime.
    A flype that flips over this inside region does not change $\widehat D$ and preserves $E$ set-wise, so we can take $D''=D$.
    A flype that flips over the outside region flips each of $D$ and $E$ along a diameter, so we can restore $E$ (if necessary) with a rotation and take $D''$ to be the result of applying this total dihedral symmetry to $D$.

    Note that $D''$ is strongly prime since it is a crossing-minimal alternating diagram (up to dihedral symmetry) for the prime tangle $T$.
    \end{proof}

    The claim can be applied to each flype in the sequence from $\widehat D$ to $\widehat D'$.
    Any flat homeomorphism $g$ in this sequence changes intermediate diagrams by planar isotopy but preserves (with labels) the two circle-components of $E$.
    Thus, $\D^2_\infty$ and $g(\D^2_\infty)$ differ only by a trivial collar and can be made to agree after a further ambient isotopy.
    The effect on $\D^2$ is therefore a (rotational) dihedral symmetry.
    At the end of this process, we have converted $D$ to $D''$ by flypes and dihedral symmetry, and we know that $\widehat D''$ and $\widehat D'$ are related by a planar isotopy.
    Now we invoke the extra structure that $f$ preserves (with labels) the two circle-components of $E$.
    This allows us to conclude that the planar isotopy preserves $\D^2_\infty$, hence it takes $D''$ onto $D'$.
    This may induce a further (rotational) dihedral symmetry.
    In the end, we have that $D$ and $D'$ are related by a sequence of flypes and dihedral symmetry, as desired.
\end{proof}

If $T$ (and $D$) are oriented, then the \emph{writhe} of $D$ is the signed count of the crossings of $D$ based on the local orientations.
Our third theorem is the analog of Theorem~\ref{thm:tait_3} and concerns the invariance of writhe among alternating diagrams for an oriented tangle.
It follows directly from Theorem~\ref{thm:alt_tang_2}, since flyping and dihedral symmetry preserve writhe.

\begin{theorem}
\label{thm:alt_tang_3}
	Any two reduced, alternating diagrams for a given oriented, prime tangle have the same writhe.
\end{theorem}

%%%%%%%%%%%%%%%%%%%%%%%%%%%%%%%%%%%%%%%%%%%%%%%%%%%%%%%%
%%%%%%%%%%%%%%%%%%%%%%%%%%%%%%%%%%%%%%%%%%%%%%%%%%%%%%%%
\subsection{Geometry of alternating tangles}
\label{sec:geometry}
%%%%%%%%%%%%%%%%%%%%%%%%%%%%%%%%%%%%%%%%%%%%%%%%%%%%%%%%
%%%%%%%%%%%%%%%%%%%%%%%%%%%%%%%%%%%%%%%%%%%%%%%%%%%%%%%%

In this section, we remark on a pair of theorems showing that alternating diagrams for tangles reflect certain aspects of the geometry of the tangles they describe.
These theorems are motivated by classical, analogous work concerning alternating diagrams for links:
The qualities of a link being prime or split are apparent in alternating diagrams.

\begin{theorem}[\cite{menasco}]
\label{thm:links_split}
    Suppose a link $L$ has an alternating diagram $D$.
    Then $L$ is prime (resp., split) if and only if $D$ is prime (resp., split).
\end{theorem}

The two theorems of interest are analogs of Theorem~\ref{thm:links_split}.
The first of these follows from the corona-closure techniques employed earlier in this section.
Recall from Section~\ref{sec:anatomy} that a split tangle is automatically composite.

\begin{theorem}
\label{thm:tangle_prime}
    Suppose a tangle $T$ has an alternating diagram $D$.
    Then $T$ is prime (resp., split) if and only if $D$ is prime (resp., split).
\end{theorem}

\begin{proof}
    If $D$ is composite (resp., split), then $T$ is composite (resp., split), as noted in Remark~\ref{rmk:diagram_notions}.
    If $T$ is composite (resp., split), then the alternating link $L$ with diagram a corona-closure $\widehat D$ of $D$ is composite (resp., split), since the decomposing sphere (resp., splitting sphere) for $T$ is automatically a decomposing sphere (resp., splitting sphere) for $L$: the inside of the sphere remains non-trivial, while the outside of the sphere gains circle-components.
    By Theorem~\ref{thm:links_split}, this implies $\widehat D$ is composite (resp., split).
    By Lemma~\ref{lem:corona}(1), this implies $D$ is composite (resp., split).
\end{proof}

\begin{remark}
\label{rmk:cromwell_tangle}
    A version of the above theorem was proved by Cromwell in the setting of tangles with two arc-components~\cite{cromwell:lonely} under a more restrictive hypothesis (strongly alternating) that involves considering closures using crossingless tangles.
    Some gaps in this approach were pointed out (and remedied) in~\cite{flapan-howards}.
    The corona-closure construction used here gives a uniform proof for tangles with any number of arc-components.

    Theorem~\ref{thm:tangle_prime} and Theorem~\ref{thm:tangle_diameters} (below) can also be proved by adapting the innermost-curve argument employed by Menasco to prove Theorem~\ref{thm:links_split}, as outlined in~\cite[Chapter~4]{Lic_97_An-introduction-to-knot-theory}.
    In the present paper, we give alternative proofs using corona-closures.
\end{remark}

\begin{figure}[ht!]
	\centering
	\includegraphics[width = .8\textwidth]{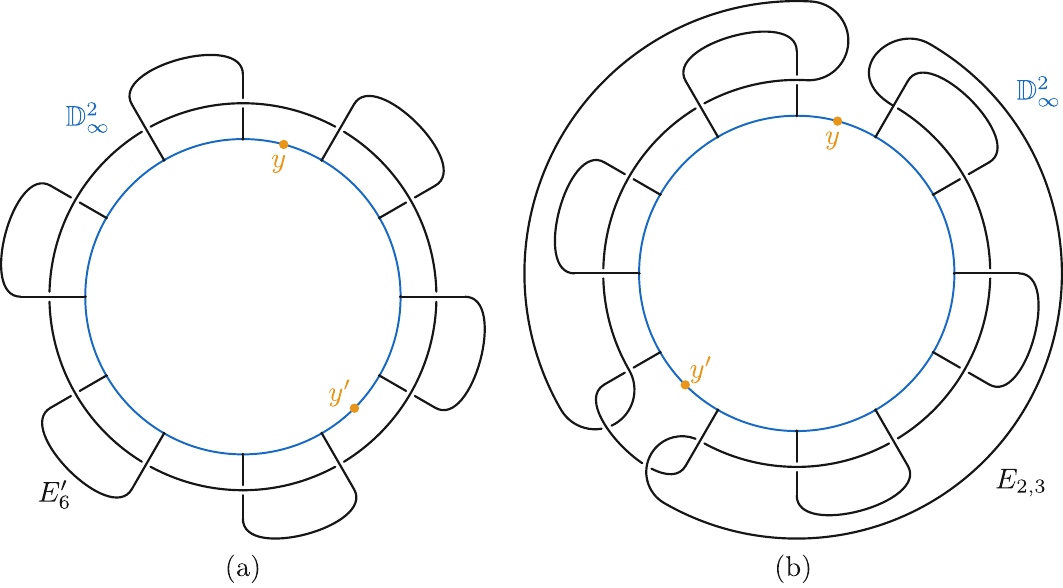}
	\caption{(a) The weak-corona diagram $E_6'$; and (b) the composite-corona diagram $E_{2,3}$.}
	\label{fig:more_coronas}
\end{figure}

For the next theorem, we will make use of variations on the corona diagram $E_n$.
Let $E'_n$ denote the alternating tangle diagram resulting from the removal of one of the circle-components from $E_n$, as shown in Figure~\ref{fig:more_coronas}(a) for $n=6$.
We will call $E_n'$ a \emph{weak-corona diagram}.
Let $E_{n_1,n_2}$ denote the alternating, diameter-composite tangle diagram whose factors are the weak-corona diagrams $E'_{n_1+1}$ and $E'_{n_2+1}$, as shown in Figure~\ref{fig:more_coronas}(b) for $n_1=2$ and $n_2=3$.
We will call $E_{n_1,n_2}$ a \emph{composite-corona diagram}.
Note that $E_n'$ and $E_{n_1,n_2}$ are visibly prime: Since $n\geq 2$ and $n_i+1\geq 2$, no pair of polygonal regions of $\D^2_\infty$ separated by the diagrams share two distinct boundary edges.

Our next result is restricted to prime, alternating tangles.

\begin{theorem}
\label{thm:tangle_diameters}
    Suppose a prime tangle $T$ has an alternating diagram $D$.
    Then $T$ is diameter-prime (resp., diameter-split) if and only if $D$ is diameter-prime (resp., diameter-split).
\end{theorem}

\begin{proof}
    First, note that one direction of each claim is straightforward (see Remark~\ref{rmk:diagram_notions}): If $D$ is diameter-composite (resp., diameter-split), then so is $T$.
    
    Conversely, first assume that $T$ is diameter-split, and let $D$ be an alternating diagram for $T$.
    Let $F$ be a vertical splitting disk for $T$, and let $y$ and $y'$ denote $F\cap\S^1$.
    Note that $y$ and $y'$ must separate the points of $T\cap\S^1$ into two sets each containing an even number of points.
    If one of these sets is empty, then the union of $F$ with a disk in $\S^2$ that is disjoint from $T$ gives a splitting sphere for $T$.
    This means $T$ is split, which means $T$ is composite, by Remark~\ref{rmk:tangle_notions}, a contradiction.

    Suppose $T$ has $n$ arc-components, and let $E'$ be either the weak-corona diagram $E_n'$ or its rotation $\rho(E_n')$ such that there is an arc $\alpha\subset\D^2_\infty$ from $y$ to $y'$ that intersects $E'$ in exactly two points.
    After potentially interchanging $T$ and $D$ with their mirrors, we can assume $\widehat D = D\cup E'$ is alternating.

    The arc $\alpha$ gives rise to a disk $F'$ whose boundary can be matched to $\partial F$ so that $F\cup F'$ is a decomposing sphere for the link $L$ whose diagram is $\widehat D$: There is at least one circle-component on either side of $F\cup F'$.
    By Theorem~\ref{thm:links_split}, $\widehat D$ is composite, with decomposing circle $\gamma$.
    Since both $D$ and $E'$ are prime, $\gamma$ cannot be moved by an ambient isotopy to be disjoint from $\S^1$.
    After an ambient isotopy, we can remove any arcs of intersection between $\gamma$ and $\D^2_\infty$ of the forms shown in Figures~\ref{fig:local}(a)--(c) (suitably adapted to a weak-corona diagram).
    After this, we can assume that $\gamma$ intersects $\S^1$ in two points.
    It follows that $\gamma = \gamma_1\cup\gamma_2$ such that $\gamma_2\subset\D^2_\infty$ intersects $E'$ in two points, and $\gamma_1\subset\D^2$ is disjoint from $D$.
    Either $\gamma_1$ is a splitting arc for $D$, as desired, or there is an ambient isotopy of $\gamma$ (through decomposing circles) moving $\gamma$ into $\D^2_\infty$, a contradiction.

    Now assume that $T$ is diameter-composite, and let $D$ be an alternating diagram for $T$.
    Note that if $D$ is diameter-split, then $D$ is diameter-composite, by a diagrammatic analog of the bottom implication in Remark~\ref{rmk:tangle_notions} (see Remark~\ref{rmk:diagram_notions}), and we are done.
    So, we assume $D$ is diameter-connected.
    Let $F$ be a vertical decomposing disk for $T$, and let $y$ and $y'$ denote $F\cap\S^1$.
    Note that $y$ and $y'$ must separate the points of $T\cap\S^1$ into two sets of odd cardinality.
    If one of these sets contains only one point, then the union of $F$ with a disk in $\S^2$ intersecting $T$ in a single point gives a decomposing sphere for $T$.
    This means $T$ is composite, a contradiction.
    Suppose these sets have size $2n_1+1$ and $2n_2+1$ for $n_1,n_2\geq 1$.

    It must be that $T$ has $n = n_1+n_2+1$ arc-components.
    Let $E = E_{n_1,n_2}$ be a composite-corona diagram such that there is an arc $\alpha\subset\D^2_\infty$ from $y$ to $y'$ that intersects $E$ in exactly one point.
    After potentially interchanging $T$ and $D$ with their mirrors, we can assume $\widehat D = D\cup E$ is alternating.

    The arc $\alpha$ gives rise to a disk $F'$ whose boundary can be matched to $\partial F$ so that $F\cup F'$ is a decomposing sphere for the link $L$ whose diagram is $\widehat D$: There is at least one circle-component on either side of $F\cup F'$.
    By Theorem~\ref{thm:links_split}, $\widehat D$ is composite, with decomposing circle $\gamma$.
    Since both $D$ and $E$ are prime, $\gamma$ cannot be moved by an ambient isotopy to be disjoint from $\S^1$.
    After an ambient isotopy, we can remove any arcs of intersection between $\gamma$ and $\D^2_\infty$ of the forms shown in Figures~\ref{fig:local}(a)--(b).
    After this, we can assume that $\gamma$ intersects $\D^2_\infty$ either as two arcs that each intersect $E$ in one point, as one arc that intersects $E$ in two points, or as one arc that intersects $E$ in one point.
    In the first two cases, $\gamma$ intersects $\D^2$ in one or two arcs that are disjoint from $D$.
    Since $D$ is diameter-connected, there is an ambient isotopy moving $\gamma$ into $\D^2_\infty$ (through decomposing circles), a contradiction.
    In the third case, $\D^2\cap\gamma$ intersects $D$ in one point, hence is a decomposing arc (since it separates $\mathbf x_{2n}$ into sets of cardinality at least three), as desired.
\end{proof}

We now prove our main theoretical result.

\begin{proof}[Proof of Theorem~\ref{thm:Tait_tangles}]
    By Theorem~\ref{thm:alt_tang_1}, we have $c(D) = c(T)$, so $c(D')\geq c(D)$.
    This yields the first claim in part (2).
    If $D'$ is reduced and alternating, then $D$ and $D'$ are flype-equivalent modulo dihedral symmetry by Theorem~\ref{thm:alt_tang_2}, establishing part (1).
    By Theorem~\ref{thm:alt_tang_1}, $c(D) = c(D') = c(T)$, establishing one direction of the second claim in part (2).
    
    Finally, assume $c(D')=c(D)$, and note that since $T$ is prime, so is $D'$, since a decomposing circle in $D'$ would immediately give a decomposing sphere for $T$; see Remark~\ref{rmk:diagram_notions}.
    By Theorem~\ref{thm:alt_tang_1}, since $D'$ is crossing-minimal, it must be alternating.
    If $D'$ is not reduced, then $c(D')>c(T)$, a contradiction.
    So, $D'$ is reduced and alternating, as desired.
\end{proof}

The following immediate corollary to Theorems~\ref{thm:Tait_tangles},~\ref{thm:tangle_prime}, and~\ref{thm:tangle_diameters} establishes that the enumeration of prime (resp., diameter-prime), connected, alternating tangles modulo isotopy rel. boundary and symmetry reduces to the enumeration of prime (resp., diameter-prime), connected, reduced, alternating tangle diagrams modulo flype-equivalence and symmetry.
Recall that $T$ is connected if and only if it is diameter-connected; see Remark~\ref{rmk:tangle_notions}.

\begin{corollary}
\label{coro:tangle_to_diag}
    Every prime, connected, alternating tangle $T$ with $c(T)=c$ modulo isotopy rel. boundary and symmetry has a unique reduced, prime, connected, alternating diagram $D$ with $c$ crossings modulo flype-equivalence and symmetry.
    Moreover, $D$ is diameter-prime if and only if $T$ is diameter-prime.
\end{corollary}

%%%%%%%%%%%%%%%%%%%%%%%%%%%%%%%%%%%%%%%%%%%%%%%%%%%%%%%%
%%%%%%%%%%%%%%%%%%%%%%%%%%%%%%%%%%%%%%%%%%%%%%%%%%%%%%%%
\section{Tangle-graphs and quadrangulations}
\label{sec:quadrangulations}
%%%%%%%%%%%%%%%%%%%%%%%%%%%%%%%%%%%%%%%%%%%%%%%%%%%%%%%%
%%%%%%%%%%%%%%%%%%%%%%%%%%%%%%%%%%%%%%%%%%%%%%%%%%%%%%%%

In this section, we rephrase the problem of enumerating prime, connected, alternating tangles in terms of disk-quadrangulations.

A \emph{tangle-graph} is an embedding $\Gamma$ of a finite graph in $\D^2$ such that
\begin{enumerate}
	\item Every vertex of $\Gamma$ in $\Int(\D^2)$ has valence four, and
	\item Every vertex of $\Gamma$ in $\S^1$ has valence one.
\end{enumerate}
Recall (Remark~\ref{rmk:n=1}) that we are assuming $n\geq 2$ and $c\geq 1$, so we have $|\S^1\cap \Gamma|\geq 4$.
It is useful to consider the graph $\Gamma\cup\S^1$, which has a ring of 3--valent vertices in place of the 1--valent vertices.
A tangle-graph will be called \emph{prime} if $\Gamma\cup\S^1$ has no 2--edge cut\footnote{Note that a prime tangle-graph cannot have a loop.}.
A connected tangle-graph will be called \emph{diameter-prime} if it cannot be disconnected by the removal of a non-leaf edge; otherwise, it will be called \emph{diameter-composite}.

\begin{lemma}
\label{lem:graph-prime}
    If $D$ is a connected, reduced, alternating tangle diagram, then $D$ is prime (resp., diameter-prime) if and only if the underlying tangle-graph $\Gamma$ is prime (resp., diameter-prime).
\end{lemma} 

\begin{proof}
    If $D$ is prime, then any simple closed curve intersecting $D$ in two non-crossing points bounds a diagram for an unknotted arc.
    Since $D$ is reduced and alternating, it is crossing-minimal, by Theorem~\ref{thm:alt_tang_1}, so the diagram for the unknotted arc must have no crossings.
    This means $D$ is strongly prime, which means $\Gamma$ is prime.
    Conversely if $D$ is composite, then the arcs of $D$ intersected by a decomposing circle give a 2--edge cut of $\Gamma\cup\S^1$.

    If $D$ is diameter-composite, then any decomposing arc for $D$ intersects $\Gamma$ in a non-leaf edge whose removal separates $\Gamma$, so $\Gamma$ is diameter-composite.
    Conversely, if $\Gamma$ is diameter-composite, then there is a diameter of $\D^2$ meeting $D$ in a single point and separating the crossings of $D$ into two non-empty sets (since $D$ is connected).
    Since $D$ is reduced and alternating, it is crossing-minimal, by Theorem~\ref{thm:alt_tang_1}, so neither side can be a diagram for an unknotted arc; hence, $D$ is diameter-composite.
\end{proof}

\begin{figure}[ht!]
	\centering
	\includegraphics[width = .8\textwidth]{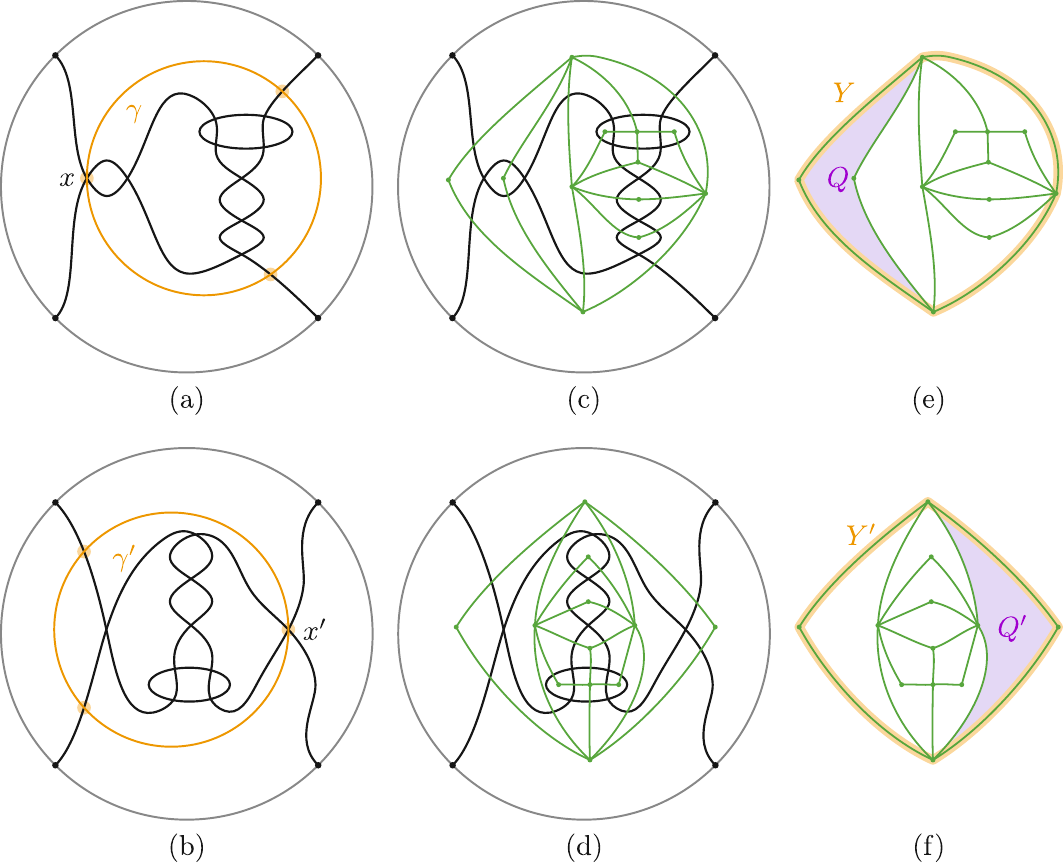}
	\caption{(a)--(b) An example of a flype along the flyping circle $\gamma$ that commutes $x$ to $x'$ and inverts the rest of the tangle-graph $\Gamma$ within $\gamma$; (c)--(d) the tangle-graph $\Gamma$ and the dual quadrangulation $G$; and (e)--(f) the effect of the flype on the region of $G$ bounded by the non-facial 4--cycle $Y$, in which the quadrilateral $Q$ is commuted from left to right and the rest of the region in $Y$ is inverted.}
	\label{fig:flypes}
\end{figure}

Two tangle-graphs are \emph{isomorphic} if there is a diffeomorphism of $\D^2$ relating them.
Flype-equivalence is defined in exactly the same way for tangle-graphs as it was for tangle diagrams; see Figures~\ref{fig:flypes}(a) and~\ref{fig:flypes}(b).
Let $\Tt$ denote the set of prime, connected tangle-graphs modulo flype-equivalence and isomorphism.
There are exactly two ways to convert a prime, connected tangle-graph into a prime, connected, reduced, alternating tangle diagram, but these two choices are mirrors of each other.
We have established the following lemma.

\begin{lemma}
\label{lem:graph_to_diagram}
    The set $\Tt$ corresponds with the set of prime, connected, reduced, alternating tangle diagrams modulo flype-equivalence and symmetry.
\end{lemma}

A \emph{disk-quadrangulation} is an embedding $G$ of a finite graph in $\D^2$ such that each face is bounded by a cycle of length four in $G$ and $\partial G = G\cap\S^1$ is a cycle in $G$.
A quadrangulation is \emph{simple} if no two edges have the same endpoints.

Given a prime, connected tangle-graph $\Gamma\subset\D^2$, define the \emph{dual quadrangulation} of $\Gamma$ to be the quadrangulation $G$ of $\D^2$ that has a face for every 4--valent vertex of $\Gamma$ and a vertex for every region of $\D^2\setminus \Gamma$\footnote{Strictly speaking, the dual quadrangulation of $\Gamma$ is a quadrangulation $G$ of a concentric sub-disk of $\D^2$, since no vertex of $G$ lies on $\S^1$; however, there is no harm in neglecting this distinction.}.
Two examples of dual quadrangulations are shown in Figure~\ref{fig:flypes}(c) and~\ref{fig:flypes}(d).
Note that $\Gamma$ is connected if and only if $\partial G$ is a cycle.

The process of passing from a connected tangle-graph to its dual quadrangulation is reversible:
Given a disk-quadrangulation $G$, each quadrilateral gives rise to a 4--valent vertex of $\Gamma$ and each edge of $\partial G$ gives rise to a leaf of $\Gamma$.

\begin{remark}
\label{rmk:numerics}
    If $\Gamma$ is connected and has $n$ arc-components and $c$ crossings, then the dual quadrangulation $G$ has $n+c+1$ vertices, $2c+n$ edges, and $c$ faces, regardless of the number $k$ of circle-components, by consideration of the Euler characteristic.
\end{remark}

If $\Gamma$ is diameter-composite, then there are two faces of $\D^2\setminus\Gamma$ that are adjacent to $\S^1$ and meet along exactly one edge $e$ that is contained in $\Int(\D^2)$.
This happens if and only if the dual quadrangulation contains an edge connecting two non-adjacent boundary vertices.
In this case, we will call the disk-quadrangulation \emph{diameter-composite}; otherwise, we will call it \emph{diameter-prime}.

A flyping circle $\gamma$ for a tangle-graph $\Gamma$ manifests as a non-facial 4--cycle $Y$ of the dual quadrangulation $G$, and the crossing $x$ incident to $\gamma$ manifests as a quadrilateral $Q$ inside $Y$ that has two edges on $Y$.
A flype along $\gamma$ manifests as a specific, simultaneous change to the portions of $\Gamma$ and $G$ inside $\gamma$ and $Y$, respectively:
The crossing $x$ and quadrilateral $Q$ are commuted to the opposite side of the inside region, and the rest of the inside region is inverted.
See Figure~\ref{fig:flypes}, where each frame in the top row is related to the corresponding frame in the bottom row by a flype along $\gamma$.

Let $\Qq$ denote the set of simple disk-quadrangulations modulo flype-equivalence and isomorphism (defined as above).

\begin{proposition}
\label{prop:tang-quad}
	The map sending a prime, connected tangle-graph to its dual quadrangulation gives a bijection $\Tt\to\Qq$.
\end{proposition}

\begin{proof}
	It is clear that two tangle-graphs are map isomorphic if and only if their dual quadrangulations are map isomorphic, and it is clear that two tangle-graphs are flype-equivalent if and only if their dual quadrangulations are flype-equivalent.
    As noted above, a tangle-graph is connected if and only if the boundary cycle of the dual quadrangulation is simple.
	
	Suppose the dual quadrangulation $G$ is not simple.
	Since a quadrangulation cannot have a loop, $G$ must have a pair of edges $e_1^*,e_2^*$ with the same endpoints.
	Let $e_1$ and $e_2$ be the edges of $\Gamma$ dual to $e_1^*,e_2^*$.
	Then, $\{e_1, e_2\}$ is a 2--edge cut for $\Gamma\cup\S^1$, so $\Gamma$ is not prime.
	Conversely, if $\Gamma$ is not prime, then it has a 2--edge cut.
    This means $G$ has a double-edge, hence is not simple.
\end{proof}

We have now reduced the problem of enumerating prime, connected, alternating tangles modulo isotopy rel. boundary and symmetry to the problem of enumerating simple disk-quadrangulations modulo flype-equivalence and isomorphism.

\begin{corollary}
\label{coro:tangle_to_quad}
    Every prime, connected, alternating tangle $T$ with $c(T)=c$ modulo isotopy rel. boundary and symmetry is represented uniquely by a simple quadrangulation $Q$ of $\D^2$ with $c$ faces modulo flype-equivalence and isomorphism.
    Moreover, $Q$ is diameter-prime if and only if $T$ is diameter-prime.
\end{corollary}

\begin{proof}
    This follows from Corollary~\ref{coro:tangle_to_diag}, Lemma~\ref{lem:graph-prime}, Lemma~\ref{lem:graph_to_diagram}, and Proposition~\ref{prop:tang-quad}.
\end{proof}

%%%%%%%%%%%%%%%%%%%%%%%%%%%%%%%%%%%%%%%%%%%%%%%%%%%%%%%%
\section{Implementation and Results}
\label{sec:implementation}
%%%%%%%%%%%%%%%%%%%%%%%%%%%%%%%%%%%%%%%%%%%%%%%%%%%%%%%%

In this section we describe the practical generation process for simple quadrangulations of $\D^2$.
First we need to revisit notions of isomorphism of maps (embedded graphs) in the sphere $\S^2$.
Two connected plane maps are \emph{isomorphic} if there is a
graph isomorphism between them that preserves the embedding
in the sphere up to orientation; i.e., reversing the orientation of
the sphere preserves the isomorphism class.
If the maps have colored edges or vertices, we also require the
isomorphism to preserve colors.

We also need a type of canonical labeling.  Let $\mathcal G$ be the set
of connected plane maps with colored vertices and edges, and let $W$ be a set
with a total order.
A \emph{signature} is a function $f:\mathcal G\to W$ such that
$G_1$ is isomorphic to $G_2$ if and only if $f(G_1)=f(G_2)$.
In our implementation (a key procedure in \texttt{plantri}), the signature is
a list of integers and our comparison is lexicographic.

Our enumeration begins with generation of \emph{disk-quadrangulations}.
As proved in the previous section, we need those which are simple and have a
simple cycle as the boundary of their outside face.
Such quadrangulations are provided by the utility \texttt{plantri} of Brinkmann and McKay~\cite{plantri,BM2007}.
Since the algorithm is not fully published, we give an overview.

First, simple quadrangulations of the sphere are generated using the algorithm of Brinkmann et al.~\cite{BGGMTW2005}.
For one vertex from each orbit of the automorphism group, we 
remove that vertex to create
a disk-quadrangulation, provided the boundary does not
have a repeated vertex.
Isomorph rejection relies on two techniques.
The canonical construction path method of McKay is used for the
quadrangulations of the sphere, ensuring that each is unique in its isomorphism
class~\cite{McKay1998}.
Next, note that two different sphere-quadrangulations can give rise to the same disk-quadrangulation---the removed vertex may have been joined to either the odd
or the even vertices on the boundary.
This ambiguity is resolved as follows.
When a vertex is removed from a sphere-quadrangulation to make a disk-quadrangulation, we color the boundary edges and also the vertices on the
boundary that were adjacent to the vertex just removed.
Call this colored map $G_1$.  Also make the colored map $G_2$ by instead
coloring the boundary vertices that were \textit{not} adjacent to the removed vertex.
Then, if $f(G_1)\le f(G_2)$ we accept~$G_1$; otherwise we reject $G_1$,
knowing that $G_2$ will be accepted instead.

Thus far we have described how \texttt{plantri} exhaustively generates disk-quadrangulations. This takes less than half a microsecond per quadrangulation
and very little space.
Next we describe how flype-equivalence is managed.
Let $\mathcal H$ be the graph whose vertices are non-isomorphic disk-quadrangulations
(assumed to have colored boundary edges), with two quadrangulations adjacent
when they are related by a single flype.
Remember that the inverse of a flype is also a flype, so this is an undirected graph.
We wish to select exactly one quadrangulation from each component of~$\mathcal H$,
and we have the tool: select the quadrangulation whose signature is least in
the component.  So, for each quadrangulation $G$ produced by \texttt{plantri},
initialize a hash table with $(G,f(G))$.  Then explore the component
of $\mathcal H$ which contains $G$ in depth-first order, inserting any quadrangulation
whose signature is new into the table and discarding any whose signature is already
in the table.
If at any stage we find a quadrangulation with signature smaller
than $f(G)$, we can abandon the DFS and discard~$G$.
If that never happens, we can write $G$ to the output.
This is efficient in practice because even if the component is
large a quadrangulation with a smaller signature is found
quickly on average.

\begin{wrapfigure}[12]{r}{0.28\textwidth}
      \centering
    \vspace{-5mm}
      \includegraphics[width=.28\textwidth]{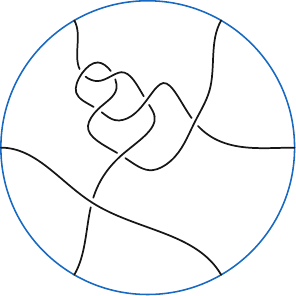}
      \caption{The tangle $T_7$.}
      \label{fig:T7}
\end{wrapfigure}
We note in passing an interesting occurrence we observed.
\begin{remark}
\label{rmk:many_flypes}
    It turns out that for each $5\leq c\leq 14$, there is a unique largest component $H_c$ of~$\Hh$.
    This component contains $2^{c-3}$ flype-equivalent disk-quadrangulations.
    
    The $c$--crossing tangle $T_c$ corresponding to $H_c$ is diameter-composite, and the diameter separates it into the one-crossing tangle and a certain rational tangle; see Figure~\ref{fig:T7} for the case of $c=7$.
    This rational tangle corresponds to the fraction $\frac{F_c}{F_{c-1}}$, where $F_c$ is the $c^\text{th}$ Fibonacci number with $F_1=F_2=1$.
    This fraction has the continued fraction expansion $[1,1,\ldots,1,2]$, where there are $c-3$ ones.

    Note that the rational tangle alone (without the extra factor making it diameter-composite) would have a smaller flype-equivalence class, since some members of the class would also be isomorphic.
    The extra arc-component breaks this symmetry.
\end{remark}

After a disk-quadrangulation is generated, we can form the corresponding tangle-graph by applying the inverse of the bijection given in Proposition~\ref{prop:tang-quad}.
We can then convert each 4--valent vertex to a crossing so that the resulting diagram is alternating.
By Corollary~\ref{coro:tangle_to_quad}, each disk-quadrangulation produced corresponds to a distinct prime, connected, alternating tangle modulo isotopy rel. boundary and symmetry.
By Theorem~\ref{thm:tangle_diameters}, the tangle is diameter-prime if and only if its alternating diagram is.
Thus, the proof of Theorem~\ref{thm:counts} is complete.

Table~\ref{tab:counts} gives counts of inequivalent prime, connected, alternating tangles as follows:
\begin{align*}
  n &= \text{the number of arc-components} \\
  k &= \text{the number of circle-components} \\
  c &= \text{the number of crossings} \\
  N_1 &= \text{the number of diameter-prime tangles} \\
  N_2 &= \text{the number of diameter-composite tangles} \\
  \text{total} &= N_1+N_2.
\end{align*}

\begin{figure}[ht!]
    \centering
    \includegraphics[width=0.875\textwidth]{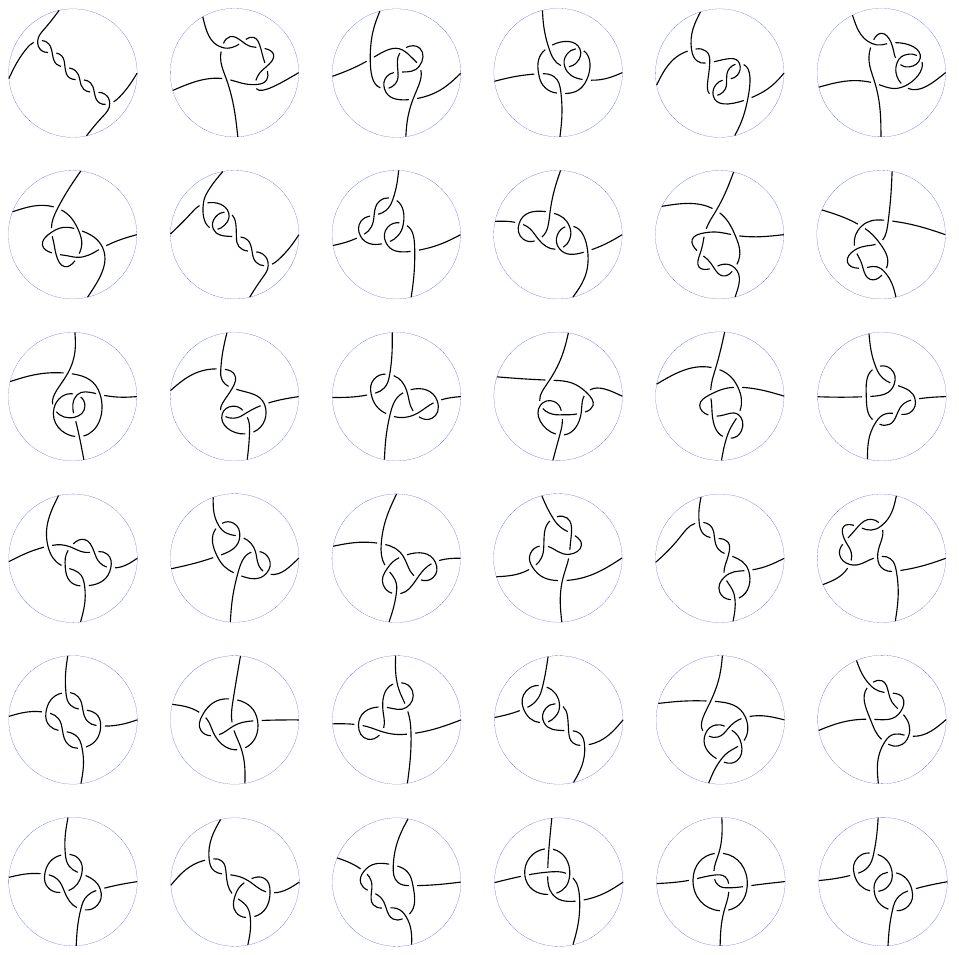}
    \caption{All prime, connected, alternating tangles with two arc-components and six crossings.}
    \label{fig:tang2_6}
\end{figure}

Values written as ``--'' are zero.
Table~\ref{tab:counts} includes the counts for $c\leq 14$, as well as for $n=2$ with $c\leq 17$, $n=3$ with $c\leq 16$, and $n=4$ with $c\leq 15$.
All the tangles listed in the table, more than 11 billion
in total, can be found in the Zenodo archive~\cite{zenodo}.
Our total counts agree with those given in Bogdanov et al., who considered $c\leq 12$ and gave diagrams for $c\leq 5$~\cite{BogMesOme_12_Enumerating-the-k-tangle-projections}.
The 36 prime, connected, alternating tangles with two arc-components and six crossings are shown
in Figure~\ref{fig:tang2_6}.
All diagrams up to 7 crossings are available in the Zenodo archive~\cite{zenodo}.

%%%%%%%%%%%%%%%%%%%%%%%%%%%%%%%%%%%%%%%%%%%%%%%%%%%%%%%%
%%%%%%%%%%%%%%%%%%%%%%%%%%%%%%%%%%%%%%%%%%%%%%%%%%%%%%%%
\bibliographystyle{amsalpha}
\bibliography{alternating_tangles.bib}
%%%%%%%%%%%%%%%%%%%%%%%%%%%%%%%%%%%%%%%%%%%%%%%%%%%%%%%%
%%%%%%%%%%%%%%%%%%%%%%%%%%%%%%%%%%%%%%%%%%%%%%%%%%%%%%%%

%\newpage
%%%%%%%%%%%%%%%%%%%%%%%%%%%%%%%%%%%%%%%%%%%%%%%%%%%%%%%%
\begin{table}[p]
\centering
\setlength{\tabcolsep}{5pt}
\begin{tabular}{rrrrrr|rrrrrr}
$n$ & $k$ & $c$ & $N_1$ & $N_2$ & total & $n$ & $k$ & $c$ & $N_1$ & $N_2$ & total \\
\hline
2 & 0 & 1 & 1 & \zero & 1 & 2 & 3 & 11 & 101 & \zero & 101 \\
2 & 0 & 2 & 1 & \zero & 1 & 2 & 3 & 12 & 707 & \zero & 707 \\
2 & 0 & 3 & 2 & \zero & 2 & 2 & 3 & 13 & 3601 & \zero & 3601 \\
2 & 0 & 4 & 4 & \zero & 4 & 2 & 3 & 14 & 20842 & \zero & 20842 \\
2 & 0 & 5 & 10 & \zero & 10 & 2 & 3 & 15 & 121413 & \zero & 121413 \\
2 & 0 & 6 & 26 & \zero & 26 & 2 & 3 & 16 & 701421 & \zero & 701421 \\
2 & 0 & 7 & 76 & \zero & 76 & 2 & 3 & 17 & 4256384 & \zero & 4256384 \\
2 & 0 & 8 & 245 & \zero & 245 & 2 & 4 & 10 & 1 & \zero & 1 \\
2 & 0 & 9 & 852 & \zero & 852 & 2 & 4 & 11 & 2 & \zero & 2 \\
2 & 0 & 10 & 3311 & \zero & 3311 & 2 & 4 & 12 & 41 & \zero & 41 \\
2 & 0 & 11 & 13667 & \zero & 13667 & 2 & 4 & 13 & 151 & \zero & 151 \\
2 & 0 & 12 & 60511 & \zero & 60511 & 2 & 4 & 14 & 1446 & \zero & 1446 \\
2 & 0 & 13 & 280128 & \zero & 280128 & 2 & 4 & 15 & 7530 & \zero & 7530 \\
2 & 0 & 14 & 1339151 & \zero & 1339151 & 2 & 4 & 16 & 53531 & \zero & 53531 \\
2 & 0 & 15 & 6582879 & \zero & 6582879 & 2 & 4 & 17 & 331908 & \zero & 331908 \\
2 & 0 & 16 & 32886971 & \zero & 32886971 & 2 & 5 & 12 & 1 & \zero & 1 \\
2 & 0 & 17 & 166994877 & \zero & 166994877 & 2 & 5 & 13 & 2 & \zero & 2 \\
2 & 1 & 4 & 1 & \zero & 1 & 2 & 5 & 14 & 57 & \zero & 57 \\
2 & 1 & 5 & 3 & \zero & 3 & 2 & 5 & 15 & 211 & \zero & 211 \\
2 & 1 & 6 & 9 & \zero & 9 & 2 & 5 & 16 & 2683 & \zero & 2683 \\
2 & 1 & 7 & 33 & \zero & 33 & 2 & 5 & 17 & 14231 & \zero & 14231 \\
2 & 1 & 8 & 110 & \zero & 110 & 2 & 6 & 14 & 1 & \zero & 1 \\
2 & 1 & 9 & 443 & \zero & 443 & 2 & 6 & 15 & 2 & \zero & 2 \\
2 & 1 & 10 & 1729 & \zero & 1729 & 2 & 6 & 16 & 78 & \zero & 78 \\
2 & 1 & 11 & 7560 & \zero & 7560 & 2 & 6 & 17 & 286 & \zero & 286 \\
2 & 1 & 12 & 34448 & \zero & 34448 & 2 & 7 & 16 & 1 & \zero & 1 \\
2 & 1 & 13 & 163884 & \zero & 163884 & 2 & 7 & 17 & 2 & \zero & 2 \\\cline{7-12}
2 & 1 & 14 & 818887 & \zero & 818887 & 3 & 0 & 2 & \zero & 1 & 1 \\
2 & 1 & 15 & 4158804 & \zero & 4158804 & 3 & 0 & 3 & 1 & 1 & 2 \\
2 & 1 & 16 & 21727691 & \zero & 21727691 & 3 & 0 & 4 & 3 & 4 & 7 \\
2 & 1 & 17 & 114731572 & \zero & 114731572 & 3 & 0 & 5 & 11 & 8 & 19 \\
2 & 2 & 6 & 1 & \zero & 1 & 3 & 0 & 6 & 42 & 26 & 68 \\
2 & 2 & 7 & 2 & \zero & 2 & 3 & 0 & 7 & 161 & 76 & 237 \\
2 & 2 & 8 & 17 & \zero & 17 & 3 & 0 & 8 & 656 & 272 & 928 \\
2 & 2 & 9 & 65 & \zero & 65 & 3 & 0 & 9 & 2776 & 1041 & 3817 \\
2 & 2 & 10 & 310 & \zero & 310 & 3 & 0 & 10 & 12321 & 4291 & 16612 \\
2 & 2 & 11 & 1477 & \zero & 1477 & 3 & 0 & 11 & 56436 & 18706 & 75142 \\
2 & 2 & 12 & 6909 & \zero & 6909 & 3 & 0 & 12 & 266791 & 84456 & 351247 \\
2 & 2 & 13 & 36097 & \zero & 36097 & 3 & 0 & 13 & 1288795 & 393286 & 1682081 \\
2 & 2 & 14 & 184457 & \zero & 184457 & 3 & 0 & 14 & 6346495 & 1875295 & 8221790 \\
2 & 2 & 15 & 1004488 & \zero & 1004488 & 3 & 0 & 15 & 31730275 & 9103836 & 40834111 \\
2 & 2 & 16 & 5500377 & \zero & 5500377 & 3 & 0 & 16 & 160667179 & 44967206 & 205634385 \\
2 & 2 & 17 & 30704047 & \zero & 30704047 & 3 & 1 & 5 & \zero & 1 & 1 \\
2 & 3 & 8 & 1 & \zero & 1 & 3 & 1 & 6 & 4 & 5 & 9 \\
2 & 3 & 9 & 2 & \zero & 2 & 3 & 1 & 7 & 17 & 21 & 38 \\
2 & 3 & 10 & 27 & \zero & 27 & 3 & 1 & 8 & 102 & 95 & 197 \\
\end{tabular}
\caption{Counts of prime, connected, alternating tangles}
\label{tab:counts}
\end{table}

\begin{table}[p]
\centering
\setlength{\tabcolsep}{5pt}
\begin{tabular}{rrrrrr|rrrrrr}
$n$ & $k$ & $c$ & $N_1$ & $N_2$ & total & $n$ & $k$ & $c$ & $N_1$ & $N_2$ & total \\
\hline
3 & 1 & 9 & 542 & 402 & 944 & 4 & 0 & 11 & 91050 & 167245 & 258295 \\
3 & 1 & 10 & 2832 & 1871 & 4703 & 4 & 0 & 12 & 468054 & 766763 & 1234817 \\
3 & 1 & 11 & 15005 & 8660 & 23665 & 4 & 0 & 13 & 2413951 & 3601338 & 6015289 \\
3 & 1 & 12 & 79578 & 42045 & 121623 & 4 & 0 & 14 & 12509991 & 17265254 & 29775245 \\
3 & 1 & 13 & 426421 & 208326 & 634747 & 4 & 0 & 15 & 65167495 & 84237466 & 149404961 \\
3 & 1 & 14 & 2303751 & 1049308 & 3353059 & 4 & 1 & 6 & \zero & 6 & 6 \\
3 & 1 & 15 & 12517758 & 5396532 & 17914290 & 4 & 1 & 7 & 2 & 45 & 47 \\
3 & 1 & 16 & 68476000 & 27991695 & 96467695 & 4 & 1 & 8 & 30 & 274 & 304 \\
3 & 2 & 7 & \zero & 1 & 1 & 4 & 1 & 9 & 244 & 1545 & 1789 \\
3 & 2 & 8 & 4 & 6 & 10 & 4 & 1 & 10 & 1797 & 8241 & 10038 \\
3 & 2 & 9 & 23 & 38 & 61 & 4 & 1 & 11 & 12151 & 43435 & 55586 \\
3 & 2 & 10 & 197 & 210 & 407 & 4 & 1 & 12 & 77674 & 227548 & 305222 \\
3 & 2 & 11 & 1226 & 1182 & 2408 & 4 & 1 & 13 & 479969 & 1196314 & 1676283 \\
3 & 2 & 12 & 7823 & 6663 & 14486 & 4 & 1 & 14 & 2894698 & 6324809 & 9219507 \\
3 & 2 & 13 & 48789 & 36414 & 85203 & 4 & 1 & 15 & 17181750 & 33651070 & 50832820 \\
3 & 2 & 14 & 297997 & 205933 & 503930 & 4 & 2 & 8 & \zero & 6 & 6 \\
3 & 2 & 15 & 1821400 & 1142588 & 2963988 & 4 & 2 & 9 & 2 & 55 & 57 \\
3 & 2 & 16 & 11014465 & 6452169 & 17466634 & 4 & 2 & 10 & 43 & 452 & 495 \\
3 & 3 & 9 & \zero & 1 & 1 & 4 & 2 & 11 & 408 & 3148 & 3556 \\
3 & 3 & 10 & 6 & 6 & 12 & 4 & 2 & 12 & 3717 & 20661 & 24378 \\
3 & 3 & 11 & 30 & 61 & 91 & 4 & 2 & 13 & 29608 & 130184 & 159792 \\
3 & 3 & 12 & 350 & 372 & 722 & 4 & 2 & 14 & 221060 & 798276 & 1019336 \\
3 & 3 & 13 & 2356 & 2720 & 5076 & 4 & 2 & 15 & 1571394 & 4825687 & 6397081 \\
3 & 3 & 14 & 17973 & 17169 & 35142 & 4 & 3 & 10 & \zero & 6 & 6 \\
3 & 3 & 15 & 125457 & 111016 & 236473 & 4 & 3 & 11 & 2 & 67 & 69 \\
3 & 3 & 16 & 865572 & 702928 & 1568500 & 4 & 3 & 12 & 60 & 680 & 740 \\
3 & 4 & 11 & \zero & 1 & 1 & 4 & 3 & 13 & 642 & 5558 & 6200 \\
3 & 4 & 12 & 7 & 8 & 15 & 4 & 3 & 14 & 6907 & 42903 & 49810 \\
3 & 4 & 13 & 39 & 91 & 130 & 4 & 3 & 15 & 61959 & 310325 & 372284 \\
3 & 4 & 14 & 584 & 620 & 1204 & 4 & 4 & 12 & \zero & 6 & 6 \\
3 & 4 & 15 & 4154 & 5468 & 9622 & 4 & 4 & 13 & 2 & 83 & 85 \\
3 & 4 & 16 & 37233 & 37795 & 75028 & 4 & 4 & 14 & 85 & 970 & 1055 \\
3 & 5 & 13 & \zero & 1 & 1 & 4 & 4 & 15 & 963 & 9130 & 10093 \\
3 & 5 & 14 & 9 & 8 & 17 & 4 & 5 & 14 & \zero & 6 & 6 \\
3 & 5 & 15 & 48 & 129 & 177 & 4 & 5 & 15 & 2 & 97 & 99 \\\cline{7-12}
3 & 5 & 16 & 922 & 965 & 1887 & 5 & 0 & 4 & \zero & 5 & 5 \\
3 & 6 & 15 & \zero & 1 & 1 & 5 & 0 & 5 & 1 & 30 & 31 \\
3 & 6 & 16 & 10 & 10 & 20 & 5 & 0 & 6 & 5 & 204 & 209 \\\cline{1-6}
4 & 0 & 3 & \zero & 2 & 2 & 5 & 0 & 7 & 39 & 1067 & 1106 \\
4 & 0 & 4 & 1 & 7 & 8 & 5 & 0 & 8 & 294 & 5415 & 5709 \\
4 & 0 & 5 & 4 & 33 & 37 & 5 & 0 & 9 & 2039 & 26269 & 28308 \\
4 & 0 & 6 & 24 & 127 & 151 & 5 & 0 & 10 & 13328 & 126729 & 140057 \\
4 & 0 & 7 & 127 & 513 & 640 & 5 & 0 & 11 & 82541 & 610721 & 693262 \\
4 & 0 & 8 & 667 & 2075 & 2742 & 5 & 0 & 12 & 492559 & 2965125 & 3457684 \\
4 & 0 & 9 & 3449 & 8686 & 12135 & 5 & 0 & 13 & 2860095 & 14528160 & 17388255 \\
4 & 0 & 10 & 17745 & 37513 & 55258 & 5 & 0 & 14 & 16305983 & 71933788 & 88239771 \\
\end{tabular}
\addtocounter{table}{-1}
\caption{Counts of prime, connected, alternating tangles (continued)}
\end{table}

\begin{table}[p]
\centering
\setlength{\tabcolsep}{5pt}
\begin{tabular}{rrrrrr|rrrrrr}
$n$ & $k$ & $c$ & $N_1$ & $N_2$ & total & $n$ & $k$ & $c$ & $N_1$ & $N_2$ & total \\
\hline
5 & 1 & 7 & \zero & 22 & 22 & 7 & 0 & 6 & \zero & 60 & 60 \\
5 & 1 & 8 & 2 & 275 & 277 & 7 & 0 & 7 & 1 & 839 & 840 \\
5 & 1 & 9 & 36 & 2190 & 2226 & 7 & 0 & 8 & 7 & 8199 & 8206 \\
5 & 1 & 10 & 473 & 15079 & 15552 & 7 & 0 & 9 & 92 & 62191 & 62283 \\
5 & 1 & 11 & 4679 & 94369 & 99048 & 7 & 0 & 10 & 1140 & 415934 & 417074 \\
5 & 1 & 12 & 40285 & 563005 & 603290 & 7 & 0 & 11 & 12370 & 2565389 & 2577759 \\
5 & 1 & 13 & 311928 & 3260729 & 3572657 & 7 & 0 & 12 & 118111 & 15075310 & 15193421 \\
5 & 1 & 14 & 2246965 & 18576455 & 20823420 & 7 & 0 & 13 & 1012431 & 85889009 & 86901440 \\
5 & 2 & 9 & \zero & 22 & 22 & 7 & 0 & 14 & 7980565 & 480029783 & 488010348 \\
5 & 2 & 10 & 2 & 353 & 355 & 7 & 1 & 9 & \zero & 612 & 612 \\
5 & 2 & 11 & 45 & 3541 & 3586 & 7 & 1 & 10 & 2 & 11178 & 11180 \\
5 & 2 & 12 & 743 & 29867 & 30610 & 7 & 1 & 11 & 64 & 124467 & 124531 \\
5 & 2 & 13 & 8829 & 224332 & 233161 & 7 & 1 & 12 & 1409 & 1103717 & 1105126 \\
5 & 2 & 14 & 89893 & 1576722 & 1666615 & 7 & 1 & 13 & 22438 & 8562582 & 8585020 \\
5 & 3 & 11 & \zero & 22 & 22 & 7 & 1 & 14 & 286421 & 61090049 & 61376470 \\
5 & 3 & 12 & 2 & 441 & 443 & 7 & 2 & 11 & \zero & 612 & 612 \\
5 & 3 & 13 & 53 & 5302 & 5355 & 7 & 2 & 12 & 2 & 14557 & 14559 \\
5 & 3 & 14 & 1125 & 52037 & 53162 & 7 & 2 & 13 & 76 & 200965 & 201041 \\
5 & 4 & 13 & \zero & 22 & 22 & 7 & 2 & 14 & 1983 & 2143563 & 2145546 \\
5 & 4 & 14 & 2 & 543 & 545 & 7 & 3 & 13 & \zero & 612 & 612 \\\cline{1-6}
6 & 0 & 5 & \zero & 16 & 16 & 7 & 3 & 14 & 2 & 18355 & 18357 \\\cline{7-12}
6 & 0 & 6 & 1 & 160 & 161 & 8 & 0 & 7 & \zero & 261 & 261 \\
6 & 0 & 7 & 6 & 1288 & 1294 & 8 & 0 & 8 & 1 & 4701 & 4702 \\
6 & 0 & 8 & 64 & 8345 & 8409 & 8 & 0 & 9 & 8 & 52807 & 52815 \\
6 & 0 & 9 & 611 & 49113 & 49724 & 8 & 0 & 10 & 133 & 456642 & 456775 \\
6 & 0 & 10 & 5361 & 272023 & 277384 & 8 & 0 & 11 & 1998 & 3401491 & 3403489 \\
6 & 0 & 11 & 42623 & 1459793 & 1502416 & 8 & 0 & 12 & 26181 & 23039389 & 23065570 \\
6 & 0 & 12 & 312623 & 7698243 & 8010866 & 8 & 0 & 13 & 293867 & 146677157 & 146971024 \\
6 & 0 & 13 & 2151049 & 40268545 & 42419594 & 8 & 0 & 14 & 2909115 & 895671650 & 898580765 \\
6 & 0 & 14 & 14107235 & 210052713 & 224159948 & 8 & 1 & 10 & \zero & 3414 & 3414 \\
6 & 1 & 8 & \zero & 119 & 119 & 8 & 1 & 11 & 2 & 72060 & 72062 \\
6 & 1 & 9 & 2 & 1749 & 1751 & 8 & 1 & 12 & 83 & 910666 & 910749 \\
6 & 1 & 10 & 51 & 16812 & 16863 & 8 & 1 & 13 & 2231 & 8953920 & 8956151 \\
6 & 1 & 11 & 841 & 132338 & 133179 & 8 & 1 & 14 & 43312 & 75847629 & 75890941 \\
6 & 1 & 12 & 10806 & 930497 & 941303 & 8 & 2 & 12 & \zero & 3414 & 3414 \\
6 & 1 & 13 & 114148 & 6110154 & 6224302 & 8 & 2 & 13 & 2 & 94161 & 94163 \\
6 & 1 & 14 & 1052265 & 38392813 & 39445078 & 8 & 2 & 14 & 98 & 1472870 & 1472968 \\
6 & 2 & 10 & \zero & 119 & 119 & 8 & 3 & 14 & \zero & 3414 & 3414 \\\cline{7-12}
6 & 2 & 11 & 2 & 2258 & 2260 & 9 & 0 & 8 & \zero & 1243 & 1243 \\
6 & 2 & 12 & 63 & 27142 & 27205 & 9 & 0 & 9 & 1 & 26752 & 26753 \\
6 & 2 & 13 & 1233 & 258428 & 259661 & 9 & 0 & 10 & 9 & 341869 & 341878 \\
6 & 2 & 14 & 19073 & 2156843 & 2175916 & 9 & 0 & 11 & 178 & 3308018 & 3308196 \\
6 & 3 & 12 & \zero & 119 & 119 & 9 & 0 & 12 & 3303 & 27155655 & 27158958 \\
6 & 3 & 13 & 2 & 2837 & 2839 & 9 & 0 & 13 & 51300 & 200166415 & 200217715 \\
6 & 3 & 14 & 72 & 40502 & 40574 & 9 & 0 & 14 & 670094 & 1371716889 & 1372386983 \\
6 & 4 & 14 & \zero & 119 & 119 & 9 & 1 & 11 & \zero & 19228 & 19228 \\
\end{tabular}
\addtocounter{table}{-1}
\caption{Counts of prime, connected, alternating tangles (continued)}
\end{table}

\begin{table}[h]
\centering
\setlength{\tabcolsep}{5pt}
\begin{tabular}{rrrrrr|rrrrrr}
$n$ & $k$ & $c$ & $N_1$ & $N_2$ & total & $n$ & $k$ & $c$ & $N_1$ & $N_2$ & total \\
\hline
9 & 1 & 12 & 2 & 466096 & 466098 & 11 & 0 & 12 & 11 & 14473264 & 14473275 \\
9 & 1 & 13 & 100 & 6584307 & 6584407 & 11 & 0 & 13 & 305 & 169415265 & 169415570 \\
9 & 1 & 14 & 3373 & 71157104 & 71160477 & 11 & 0 & 14 & 7933 & 1646170373 & 1646178306 \\
9 & 2 & 13 & \zero & 19228 & 19228 & 11 & 1 & 13 & \zero & 650325 & 650325 \\
9 & 2 & 14 & 2 & 610916 & 610918 & 11 & 1 & 14 & 2 & 19673444 & 19673446 \\\cline{1-6}\cline{7-12}
10 & 0 & 9 & \zero & 6257 & 6257 & 12 & 0 & 11 & \zero & 175760 & 175760 \\
10 & 0 & 10 & 1 & 155592 & 155593 & 12 & 0 & 12 & 1 & 5464660 & 5464661 \\
10 & 0 & 11 & 10 & 2221534 & 2221544 & 12 & 0 & 13 & 12 & 94498320 & 94498332 \\
10 & 0 & 12 & 239 & 23758271 & 23758510 & 12 & 0 & 14 & 390 & 1201387390 & 1201387780 \\
10 & 0 & 13 & 5219 & 212946156 & 212951375 & 12 & 1 & 14 & \zero & 3857256 & 3857256 \\\cline{7-12}
10 & 0 & 14 & 94897 & 1696596887 & 1696691784 & 13 & 0 & 12 & \zero & 963900 & 963900 \\
10 & 1 & 12 & \zero & 111111 & 111111 & 13 & 0 & 13 & 1 & 32885537 & 32885538 \\
10 & 1 & 13 & 2 & 3024387 & 3024389 & 13 & 0 & 14 & 13 & 618158135 & 618158148 \\\cline{7-12}
10 & 1 & 14 & 123 & 47229498 & 47229621 & 14 & 0 & 13 & \zero & 5374400 & 5374400 \\
10 & 2 & 14 & \zero & 111111 & 111111 & 14 & 0 & 14 & 1 & 199535742 & 199535743 \\\cline{1-6}\cline{7-12}
11 & 0 & 10 & \zero & 32721 & 32721 & 15 & 0 & 14 & \zero & 30385256 & 30385256 \\
11 & 0 & 11 & 1 & 916594 & 916595 &  &  &  &  &  &  \\
\end{tabular}
\addtocounter{table}{-1}
\caption{Counts of prime, connected, alternating tangles (continued)}
\end{table}
%%%%%%%%%%%%%%%%%%%%%%%%%%%%%%%%%%%%%%%%%%%%%%

\end{document}